%% file: Localization1.tex
\documentclass{amsart}

\usepackage{amsmath,amssymb,amsfonts,amsthm,amscd,latexsym,euscript}
\usepackage[all]{xy}

\input{symbols.tex}

\begin{document}

\SelectTips{cm}{10}

%
%

\title[Equivariant localization of diagrams]{Localization with respect to a class of maps I --
                                            Equivariant localization of diagrams of spaces}
\author{Boris Chorny}
\address{Einstein Institute of Mathematics, Edmond Safra Campus, Givat
Ram, The Hebrew University of Jerusalem, Jerusalem, 91904, Israel}
\curraddr{Department of Mathematics, Middlesex College, The University of Western Ontario, London, Ontario N6A 5B7, Canada}

\email{bchorny2@uwo.ca}

\thanks{During the preparation of this paper the author was a fellow of
        Marie Curie Training Site hosted by Centre de Recerca Matem\`atica
       (Barcelona), grant no. HPMT-CT-2000-00075 of the European Commission.}

\subjclass{Primary 55U35; Secondary 55P91, 18G55} \keywords{model category,
localization, equivariant homotopy}
\date{\today}
\dedicatory{}
\commby{}

\hfuzz=2.5pt

\begin{abstract}
Homotopical localizations with respect to a set of maps are known
to exist in cofibrantly generated model categories (satisfying
additional assumptions) \cite{Bous:factor, Farjoun-book,
Hirschhorn, Smith}. In this paper we expand the existing
framework, so that it will apply to not necessarily cofibrantly
generated model categories and, more important, will allow for a
localization with respect to a class of maps (satisfying some
restrictive conditions).

We illustrate our technique by applying it to the equivariant
model category of diagrams of spaces \cite{Farjoun}. This model
category is not cofibrantly generated \cite{Chorny}. We give
conditions on a class of maps which ensure the existence of the
localization functor; these conditions are satisfied by any set of
maps and by the class of maps which induces ordinary localizations
on the generalized fixed-points sets.
\end{abstract}

\maketitle


\section*{Introduction}

Homotopy idempotent constructions, or homotopical localizations,
play an important role in algebraic topology and algebraic
geometry. Homotopical localization is a functor $L$ in a model
category that carries weak equivalences into weak equivalences and
is equipped with the natural transformation $Id\rarrow L$, that
induces weak equivalences $LX\simeq LLX$ for all $X$.

The idea behind such construction is to ``forget information'' in
a consistent, functorial way. The amount of the information
``forgotten'' by $L$ is measured precisely by the class of maps
$\cal S_L$ which are turned into weak equivalences by $L$.

In practice, we usually know what kind of information we would like
to discard, i.e. $\cal S_L$, and we are looking for the functor $L$.
E.g., $\cal S_L = \{$homological equivalences of spaces$\}$;
and $L$ is Bousfield's localization functor \cite{Bous:homology}.

If we are able to encode the ``informative content'' of $\cal S_L$ by just a 
set of maps $S\subset \cal S_L$, and if the underlying model category is
cofibrantly generated (plus some other conditions), then the construction of
the localization functor $L$ is given by the classical framework established by
A.K~Bousfield, E.~Dror Farjoun, P.S.~Hirschhorn, J.~Smith \cite{Bous:factor,
Farjoun-book, Hirschhorn, Smith}.

But if we need to invert a proper class of maps or if we happen to work in a
non-cofibrantly generated model category, then the standard framework does not
guarantee the existence of the localization functor. Such situations are not
rare. Recently several important model categories were shown to be
non-cofibrantly generated \cite{AHRT, Chorny, ChristHov, Isaksen}. Also the
existence of the localization functors with respect to the class of
cohomological equivalences of spaces is a long-standing open problem.

In this work we develop a new approach that allows for construction of
homotopical localizations with respect to a class of maps (satisfying some
conditions) in a not necessarily cofibrantly generated model category. As an
application, we give an example of a non-cofibrantly generated model category (taken from \cite{Chorny}) and a class of maps which fits into our framework but is not covered by the
classical localization frameworks. In our example the model structure is generated by classes of cofibrations and trivial cofibrations which carry an additional structure sufficient to guarantee the applicability of a version of Quillen's small object argument. The set of conditions on a class of maps (in an abstract category) that allows one to apply the generalized version of the argument is given the name \emph{instrumentation}.

The main example of this paper involves homotopical localizations in the
equivariant homotopy theory. Classical Bredon's homotopy theory of spaces with a group action was generalized to the equivariant homotopy theory over an arbitrary small category $D$ by E.~Dror Farjoun and A.~Zabrodsky in the mid-80's \cite{DZ, DF, Farjoun}. 

The question of the existence of equivariant localization functors was
considered previously by V.~Halperin \cite{Halperin}. Only \emph{strong}
equivariant localizations were constructed in that work, which are not
localizations in any familiar model category on diagrams of spaces.

\subsection*{Organization of the paper} 
In the preliminary section we recall basic definition and results about the homotopy theory of diagrams of spaces. We also make explicit the connection between two notions of orbits: Dwyer-Kan \cite{DK} and Dror Farjoun-Zabrodsky \cite{DZ, Farjoun}. Namely the orbits in the later works are also orbits from the first work. 

In Section~2 we develop a formalism necessary to describe the additional structure carried by a class of maps in order to satisfy the conditions of Quillen's small object argument. We introduce here the central notion of \emph{instrumentation} of a class of maps. In Section~3 we illustrate this notion by constructing instrumentations on the classes of generating cofibrations and trivial cofibrations in the category of diagrams of spaces equipped with the equivariant model structure of \cite{Farjoun}.

Section~4 is devoted to the proof of the generalized Quillen small object argument, which also finishes the proof that the factorizations in the model category of \cite{Farjoun} are functorial.

Before turning to the construction of localizations in the equivariant model category of $D$-shaped diagrams of spaces we prove in Section~5 that this model category is proper. We also answer affirmatively the question posted by E.~Dror Farjoun in \cite[2.3]{Farjoun} and show that in the category of diagrams of simplicial sets every object is cofibrant.

In Section~6 we list the properties required from a class of maps in the equivariant category of diagrams in order to ensure the existence of the localization functor with respect to that class of maps. Afterwards we construct these localization functors and prove their basic properties.

Finally, in Section~7 we apply the technique of Section~6 in order to construct the fixed-pointwise localization functor with respect to a map of spaces (for the case of $G$-spaces, where $G$ is a compact Lie group, such localizations were constructed by J.~P.~May, J.~McClure and G.~Triantafillou \cite{MMT} with respect to the ordinary non-equivariant homology theory).

In Appendix~A we discuss the notion of contractible objects in a general model category. Their properties were useful in the proof of the properties of localization functors in Section~6.

\medskip
\paragraph{\emph{Acknowledgements.}} I would like to thank Emmanuel Dror Farjoun for
his support and many helpful ideas. I am grateful to Carles Casacuberta  who
has patiently read and helped me to improve the early version of this paper.

This paper was significantly influenced by the fundamental treatise of
P.~S.~Hir\-sch\-horn \cite{Hirschhorn}. I am obliged to the author who made
available through the internet the early versions of his manuscript.

\section{Preliminaries on the diagrams of spaces}
In this section we review the equivariant homotopy theory of diagrams of
spaces and establish basic notation. The only novelty introduced here is that
an orbit diagram in the sense of Dror Farjoun--Zabrodsky \cite{DF,Farjoun, DZ}
is also an orbit in the sense of Dwyer--Kan \cite{DK}. The topological (versus
simplicial) version of the main proposition appeared previously in
\cite{Chorny1}.

In this paper the \emph{category of spaces} ${\cal S}$ is the category of
simplicial sets or compactly generated topological spaces with the standard
model structure. Most of the results are true in the category of pointed
spaces. If $D$ is a small category enriched over $\cal S$, then the
\emph{category of ($D$-shaped) diagrams of spaces} ${\cal S}^D$ is the
category of continuous functors from $D$ to $\cal S$ with natural
transformations as morphisms. ${\cal S}^D$ is a simplicial category and we
denote by $\hom(\cdot, \cdot)$ the simplicial function complex; $\hom_D(\cdot,
\cdot)$ will denote the set of morphisms in $\cal S^D$. The two functors are related as follows:
\[
\hom(\dgrm X, \dgrm Y)_n = \hom_D(\dgrm X\otimes \Delta^n, \dgrm Y).
\]

There are several well-known model structures on the category of diagrams of
spaces. One of the most widely used is the Bousfield--Kan model category, in
which weak equivalences and fibrations are objectwise and cofibrations are
obtained by the left lifting property with respect to trivial fibrations.
Another example (not used in this paper) is A.~Heller's model category, in which weak
equivalences and cofibrations are objectwise and fibrations are obtained by
the right lifting property with respect to trivial cofibrations. These model
categories are cofibrantly generated. In this article we work mostly with the equivariant model structure on the category of diagrams of spaces constructed by E.~Dror Farjoun \cite{Farjoun} and described in Definition~\ref{model-structure}.

Recall from \cite{DK} that a set $\{O_e\}_{e\in E}$ of objects of a category
\cat M, enriched over simplicial sets, is said to be a set of orbits for \cat M
if the following axioms hold:
\begin{itemize}
\item[Q0:] \cat M is closed under arbitrary direct limits.
\item[Q1:]  For every $e\in E$, the functor $\hom(O_e,\cdot)\colon \cat M \rarrow \cal
Sets^{\Delta^{\text{op}}}$ commutes, up to homotopy, with the pushouts of the
following form:
\[
\xymatrix{
O_{e'}\otimes K \ar@{^{(}->}[d] \ar[r]     &       X_a\ar[d]\\
O_{e'}\otimes L \ar[r]                     &       X_{a+1}
               }
\]
where $K\hookrightarrow L$ is an inclusion of finite simplicial sets.
\item[Q2:] For every $e\in E$, the functor $\hom(O_e,\cdot)\colon \cat M \rarrow \cal
Sets^{\Delta^{\text{op}}}$ commutes, up to homotopy, with transfinite
compositions of maps $X_a\hookrightarrow X_{a+1}$ as in Q1.
\item[Q3:] There is a limit ordinal $\kappa$ such that, for every $e\in E$, the
functor $\hom(O_e, \cdot)$ strictly commutes with $\kappa$-transfinite
compositions of maps $X_a\hookrightarrow X_{a+1}$ as in Q1.
\end{itemize}

According to \cite{DF,Farjoun,DZ}, a $D$-diagram $\dgrm T$ is an \emph{orbit}
if $\colim_D \dgrm T = \ast$; $\cat O_D$ denotes the full subcategory of
orbits. Often $\cat O_D$ is a large subcategory of $\cal S^D$. The objective
of this section is to show that the orbits in this sense satisfy the axioms
Q0--Q3. Nevertheless, the collection of orbits is not a set of orbits for
$\cal S^D$, since it may be a proper class rather than a set. Any subset of
$\obj{\cat O_D}$ is a set of orbits for $\cal S^D$ and defines a model
category structure; see \cite{DK}.

Axiom Q0 is obvious in the category of diagrams of spaces. Axiom Q2 was
verified in \cite{Farjoun}. Axioms Q1 and Q3 were left to the reader in \cite{Farjoun}. We decided to include them into current exposition. Axiom Q3 will be proved in Proposition~\ref{small:domains} below. The following proposition verifies Axiom Q1.

\begin{proposition}\label{Q1}
If \dgrm T is an orbit, $K\hookrightarrow L$ is an inclusion of (finite)
simplicial sets and the square
\begin{equation} \label{pushoutQ1}
\xymatrix{
\dgrm T\otimes K \ar@{^{(}->}[d] \ar[r]     &       \dgrm X_a\ar[d]\\
\dgrm T\otimes L \ar[r]                     &       \dgrm X_{a+1}
               }
\end{equation}
is a pushout diagram in $\cal S^D$, then for any other orbit $\dgrm T'$
the commutative square
\[
\xymatrix{
 \hom(\dgrm T', \dgrm T\otimes K) \ar@{^{(}->}[d] \ar[r]     &       \hom(\dgrm T', \dgrm X_a)\ar[d]\\
 \hom(\dgrm T', \dgrm T\otimes L) \ar[r]                     &       \hom(\dgrm T', \dgrm X_{a+1})
               }
\]
is a pushout diagram, up to homotopy, in the category of simplicial sets, i.e.,
the natural map $\hom(\dgrm T', \dgrm X_a) \coprod_{\hom(\dgrm T', \dgrm
T\otimes K)}\hom(\dgrm T', \dgrm T\otimes L) \rarrow \hom(\dgrm T', \dgrm
X_{a+1})$ is a weak equivalence of simplicial sets.
\end{proposition}
\begin{proof}
The argument for $\cal S = \Top$ differs from the case $\cal S = \cal Sets^{\Delta^{\text{op}}}$. We will treat first the case of the simplicial sets.

Let us prove first the following private case: for any orbit $\dgrm T'$, the functor $\hom(\dgrm T', \cdot)$ strictly commutes with the pushouts of the form 
\begin{equation}\label{pushoutQ1'}
\xymatrix{
\dgrm T\otimes K \ar@{^{(}->}[d] \ar[r]     &       \dgrm X_a\ar@{^{(}->}[d]\\
\dgrm T\otimes L \ar@{^{(}->}[r]            &       \dgrm X_{a+1},
               }
\end{equation}
where $\hookrightarrow$ are cofibration in the sense of Definition~\ref{model-structure} bellow, or transfinite compositions of the maps as in Q1 above. Note that the restriction applies only on the lower horizontal map; the right vertical map is a cofibration also in (\ref{pushoutQ1}).

It will suffice to show that in each dimension $n\geq 0$ the commutative square
\[
\xymatrix{
 \hom_D(\dgrm T'\otimes \Delta^n, \dgrm T\otimes K) \ar[d] \ar[r]     &       \hom_D(\dgrm T'\otimes \Delta^n, \dgrm X_a)\ar[d]\\
 \hom_D(\dgrm T'\otimes \Delta^n, \dgrm T\otimes L) \ar[r]                     &       \hom_D(\dgrm T'\otimes \Delta^n, \dgrm X_{a+1})
               }
\]
is a pushout in the category of sets. The functor of tensoring with a
simplicial set $W$ is equal to the product, in the category of diagrams, with
the constant diagram containing $W$ in each entry. Hence, it commutes with
$\hom_D(\dgrm T'\otimes \Delta^n,\;\cdot\;)$. Additionally, we recall that the
tensor $\cdot\otimes \Delta^n$ is the left adjoint of the cotensor functor
$(\;\cdot\;)^{\Delta^n}$, therefore the commutative square above becomes
\[
\xymatrix{
 \hom_D(\dgrm T'\otimes \Delta^n, \dgrm T)\times \hom_D(\dgrm T', K^{\Delta^n}) \ar[d] \ar[r]     &       \hom_D(\dgrm T'\otimes{\Delta^n}, \dgrm X_a)\ar[d]\\
 \hom_D(\dgrm T'\otimes \Delta^n, \dgrm T)\times \hom_D(\dgrm T', L^{\Delta^n}) \ar[r]                     &       \hom_D(\dgrm T'\otimes{\Delta^n}, \dgrm X_{a+1})
               }
\]
where the simplicial sets $K^{\Delta^n}$ and $L^{\Delta^n}$ are thought of as
constant diagrams.

Let $C$ be the set of connected components of the nerve of $D$. Since $\dgrm
T'$ is an orbit, any map from $\dgrm T'$ to a constant diagram with simplicial
set $W$ in each entry is determined by the image of $\colim_D \dgrm T' =
\Delta^0$ in $\colim_D W= \coprod_C W$, but $\Delta^0$ can hit only one
component, determined by $\dgrm T'$, i.e., $\hom_D(\dgrm T', W) = W_0$ -- the set of 0-simplices of $W$. We conclude that
\[
\hom_D(\dgrm T', K^{\Delta^n}) = (K^{\Delta^n})_0 = K_n,\quad 
\hom_D(\dgrm T', L^{\Delta^n}) = (L^{\Delta^n})_0 = L_n. 
\]

Additionally, note that $\colim_D \dgrm X_{a+1} = \colim_D \dgrm X_{a}
\coprod_{K}L$, and, on the level of the $n$-simplices, $(\colim_D \dgrm
X_{a+1})_n =(\colim_D \dgrm X_{a})_n\coprod_{K_n} L_n = (\colim_D \dgrm
X_{a})_n\coprod (L_n\setminus K_n)$. Then the set of maps from $\dgrm T'\otimes
{\Delta^n}$ to $\dgrm X_{a+1}$ may be decomposed into the disjoint union of sets parameterized by the $(\colim_D\dgrm X_{a+1})_n$ as follows:
\begin{multline*}
\hom_D(\dgrm T'\otimes{\Delta^n}, \dgrm X_{a+1}) =\\
\hom_D(\dgrm T'\otimes{\Delta^n}, \colim_D \dgrm X_{a}\times_{\underset{D}{\colim} \dgrm X_{a+1}}\dgrm X_{a+1})\coprod 
\left(\coprod_{x\in L_n\setminus K_n} \hom_D^x(\dgrm T'\otimes {\Delta^n}, \dgrm P_x)\right ),\\
\end{multline*}
where $\dgrm P_x$ is the pullback of the map $\dgrm X_{a+1}\rarrow \colim_D \dgrm X_{a+1}$ over the $n$-simplex $x\in L_n\setminus K_n$, i.e., the map $x\colon \Delta^n\rarrow \colim_D \dgrm X_{a+1}$; and $\hom_D^x$ is the set of all equivariant maps which induce the map $x$ on the colimits.

The restriction imposed on the commutative square (\ref{pushoutQ1'}) and \cite[Lemma 2.1]{Farjoun} imply that the following commutative squares are pullbacks:
\[
\xymatrix{
\dgrm T\otimes L \ar@{^{(}->}[r] 
                 \ar[d]  &  \dgrm X_{a+1}\ar[d]    & & \dgrm X_{a}\ar@{^{(}->}[r]
                                                                  \ar[d] &\dgrm X_{a+1}\ar[d]\\
          L      \ar@{^{(}->}[r]  &  \colim_D \dgrm X_{a+1},
                                                   & & \colim_D \dgrm  X_{a}\ar@{^{(}->}[r] &
                                                                       \colim_D \dgrm  X_{a+1}.
               }
\]

Therefore, $\dgrm P_x = \dgrm T\times \Delta^n$ and $\colim_D \dgrm X_{a}\times_{\underset{D}{\colim} \dgrm X_{a+1}}\dgrm X_{a+1} = \dgrm X_a$. Finally we can conclude:
\[
\hom_D(\dgrm T'\otimes{\Delta^n}, \dgrm X_{a+1}) = \hom_D(\dgrm T'\otimes{\Delta^n}, \dgrm X_{a})\coprod \left(\hom_D(\dgrm T'\otimes {\Delta^n}, T)\times (L_n\setminus K_n)\right)
\]

Assembling all the information obtained so far, we rewrite the initial
commutative square as

\[
\xymatrix{
 \hom_D(\dgrm T'\otimes \Delta^n, \dgrm T)\times K_n \ar@{^{(}->}[d] \ar[r]     &       \hom_D(\dgrm T'\otimes \Delta^n, \dgrm X_a)\ar[d]\\
 \hom_D(\dgrm T'\otimes \Delta^n, \dgrm T)\times L_n \ar[r]                     &       \scriptstyle{\hom_D(\dgrm T'\otimes{\Delta^n}, \dgrm X_{a})\coprod \hom_D(\dgrm T'\otimes{\Delta^n}, T)\times (L_n\setminus K_n).}
}
\]
We need to show that this square is a pushout of sets. It is implied by the
following decomposition of this square into the disjoint union of two pushout
squares:
\[
\xymatrix{
 \hom_D(\dgrm T'\otimes \Delta^n, \dgrm T)\times K_n \ar@{=}[d] \ar[r]    &       \hom_D(\dgrm T'\otimes \Delta^n, \dgrm X_a)\ar@{=}[d]\\
 \hom_D(\dgrm T'\otimes \Delta^n, \dgrm T)\times K_n \ar[r]               &       \hom_D(\dgrm T'\otimes{\Delta^n}, \dgrm X_{a})
               }
\]
and
\[
\xymatrix{
 \emptyset \ar[d] \ar@{=}[r]     &       \emptyset\ar[d]\\
 \hom_D(\dgrm T'\otimes \Delta^n, \dgrm T)\times (L_n\setminus K_n) \ar@{=}[r]       &       \hom_D(\dgrm T'\otimes{\Delta^n}, T)\times (L_n\setminus K_n).
               }
\]
Therefore, the functor $\hom(\dgrm T, \;\cdot\;)$ commutes with pushouts of the
form (\ref{pushoutQ1'}).

For the general case, factor the map $\dgrm T\otimes L\rarrow \dgrm X_{a+1}$ into a cofibration followed by a trivial fibration (this is legitimate from the point of view of \cite{Farjoun} since Q1 is only required for the construction of the second factorization; see \cite[2.2]{DK})
\[
\xymatrix{
\dgrm T\otimes K \ar@{^{(}->}[d] \ar[r] &  \dgrm X_{a}'\ar@{->>}[r]^\sim 
                                                       \ar[d]  &    \dgrm X_a\ar@{^{(}->}[d]\\
\dgrm T\otimes L \ar@{^{(}->}[r]        &  \dgrm X_{a+1}'\ar@{->>}[r]^\sim
                                                               &    \dgrm X_{a+1},
               }
\]
and take $\dgrm X_a'= \dgrm X_{a+1}'\times_{\dgrm X_{a+1}}\dgrm X_{a}$. Then the left square is a pushout: $\dgrm X_{a+1}' = \dgrm X_{a}'\coprod_{\dgrm T\otimes K}\dgrm T\otimes L$ by the topos theory result \cite[IV.7.2]{MacMoer} applied to the topoi $\cal S^D\overcat\dgrm X_{a+1}'$ and $\cal S^D\overcat \dgrm X_{a+1}$ (recall that $\cal S = \cal Sets^{\Delta^\op}$ is a topos), since it may be obtained by the base change from the outer square along the map $\dgrm X_{a+1}'\rarrow \dgrm X_{a+1}$.

Then the left square is of the form (\ref{pushoutQ1'}) and in the diagram
\[
\xymatrix{
\hom(\dgrm T',\dgrm T\otimes K) \ar@{^{(}->}[d] \ar[r] &  \hom(\dgrm T',\dgrm X_{a}')\ar@{->>}[r]^\sim 
                                                       \ar[d]  &    \hom(\dgrm T',\dgrm X_a)\ar[d]\\
\hom(\dgrm T',\dgrm T\otimes L) \ar[r]        &  \hom(\dgrm T',\dgrm X_{a+1}')\ar@{->>}[r]^\sim
                                                               &    \hom(\dgrm T',\dgrm X_{a+1}),
               }
\]
the left square is a pushout and the outer square is a pushout up to homotopy, as required. 

In the case $\cal S = \Top$, there is the topological function space
$\hom^{\Top} \colon (\Top^D)^{\text{op}}\times \Top^D \rarrow \Top$. In order
to pass from the topological to simplicial function space one applies the
singular functor. It was proven in \cite[Proposition~3.1]{Chorny1} that the
functor $\hom^{\Top}(\dgrm T', \;\cdot\;)$ commutes with the pushouts of the
form (\ref{pushoutQ1}). This means that the commutative square
\[
\xymatrix{
 \hom^{\Top}(\dgrm T', \dgrm T\otimes K) \ar@{^{(}->}[d]_i \ar[r]   &       \hom^{\Top}(\dgrm T', \dgrm X_a)\ar[d]\\
 \hom^{\Top}(\dgrm T', \dgrm T\otimes L) \ar[r]                     &       \hom^{\Top}(\dgrm T', \dgrm X_{a+1})
               }
\]
is a pushout in the category of compactly generated topological spaces.

We need to show that the application of the singular functor preserves this
pushout up to homotopy, though it is not true in general. It follows from
\cite[13.5.5]{Hirschhorn}, since the map $i$ is a cofibration ($i$ is equal to
$\hom^{\Top}(\dgrm T', \dgrm T)\times |K| \hookrightarrow \hom^{\Top}(\dgrm
T', \dgrm T)\times |L|$).
\end{proof}

Although the collection of orbits is not, in general, a set of orbits for $\cal
S^D$, it determines a model structure given by the following
\begin{definition}\label{model-structure}
We say that a model structure on the category of the $D$-shaped diagrams of
spaces is \emph{generated by a collection $\cal O$ of orbits} if a map $f\colon\dgrm X
\rarrow \dgrm Y$ is
\begin{itemize}
\item
a \emph{weak equivalence} if and only if the induced map
\[
\hom(\dgrm T, f) \colon \hom(\dgrm T, \dgrm X) \rarrow \hom(\dgrm T, \dgrm Y)
\]
is a weak equivalence of simplicial sets for any orbit $\dgrm T\in \cal O$;
\item
a \emph{fibration} if and only if the induced map
\[
\hom(\dgrm T, f) \colon \hom(\dgrm T, \dgrm X) \rarrow \hom(\dgrm T, \dgrm Y)
\]
is a fibration of simplicial sets for any orbit $\dgrm T\in \cal O$;
\item
a \emph{cofibration} if and only if it has the left lifting property with
respect to trivial fibrations.

We say that the model structure is \emph{equivariant} if $\cal O = \obj{\cat O_D}$.
\end{itemize}
\end{definition}

\begin{example}
If $D=G$ is a group, then a $G$-orbit is just a homogeneous space $G/H$ for
some subgroup $H<G$ (hence the name -- orbit). In this case the collection of
orbits forms a set $\{G/H \;|\; H<G\}$. Weak equivalences from
Definition~\ref{model-structure} coincide with the classical $G$-equivariant
homotopy equivalences introduced by G.~Bredon in \cite{Bredon}. The equivariant model
category generated by the set of orbits was constructed in \cite{DK}.
\end{example}
\begin{example}
If $D = (\bullet \rarrow \bullet)$ is the category with two objects and one
non-identity morphism, then for any space $X\in \cal S$ there corresponds the
orbit $T_X = (X \rarrow \ast)$. The full subcategory of orbits in this case is
equivalent to the category of spaces \cal S.
\end{example}

Although the full homotopical information in the equivariant model category is encoded, usually, by the class of orbits, for a fixed diagram there exists a set of orbits, which captures its homotopical information.

\begin{definition}
The \emph{category of orbits} of \dgrm X is a full small subcategory of $\cat
O_D$ generated by the set of orbits
\[
\{\dgrm T_x = \ast \times_{\colim_D \dgrm X} \dgrm X \;|\; x\colon\ast\rarrow
\colim_D \dgrm X\},\; \text{for $\cal S = \Top$};
\]
or
\[
\{\dgrm T_x = \Delta^0 \times_{\colim_D \dgrm X^{\Delta^n}} \dgrm X^{\Delta^n}
\;|\; x\colon \Delta^0 \rarrow \colim_D \dgrm X^{\Delta^n}, n\geq 0\},\;
\text{for $\cal S = \Sets^{\Delta^{\op}}$}.
\]
The category of orbits of \dgrm X is denoted by $\cat O_{\dgrm X}$.
\end{definition}

Let $\cat E \subset \cat O_D$ be a small full subcategory. Then $\obj{\cat E}$
forms a set of orbits for $\cal S^D$. For any diagram \dgrm X, we can form a
diagram $\dgrm X^{\cat E}$ of spaces over the category $\cat E^\op$ with
$\dgrm X^{\cat E}(\dgrm E) = \hom(\dgrm E, \dgrm X), \; \dgrm E\in \obj{\cat
E}$. We will call it the diagram of orbit-points of \dgrm X, since it
generalizes the notion of the diagram of fixed-points of a space with a group
action. According to \cite{DK}, there exists a model category on $\cal S^D$
with a map $f$ being a weak equivalence or fibration if and only if $f^{\cat E}$
is a weak equivalence or fibration in the Bousfield--Kan model category on
$\cal S^{\cat E^\op}$. Moreover, the functor $(\;\cdot\;)^{\cat E}$ has a left
adjoint and this is a Quillen equivalence of the model categories.

Fix a diagram \dgrm X, and assume that \dgrm X is of orbit type $\cat E$,
i.e., $\cat O_{\dgrm X}\subset \cat E$. We will use the result
\cite[Lemma~3.7]{DF}, which shows that $\dgrm X^{\cat E}$ is $\cat E^\op$-free.
This means, essentially, that the orbit-points functor preserves cofibrant
diagrams of orbit type $\cat E$, which is not typical for a right adjoint.

\section{Instrumented classes of maps}
In this section we introduce several notions from Category Theory. Some of
them are well-known, others are new. The main new concept presented here
(Definition~\ref{instr}) is the \emph{instrumented class of maps}. This notion
will allow us to generalize Quillen's small object argument in Section
\ref{Quillen}. Any set of morphisms with small domains (see below) in any
category may be thought of as an instrumented class. A non-trivial example of
an instrumented (proper) class of maps in the category of diagrams of spaces
is given in Section~\ref{Diagrams}.
\par
Throughout this section, let \cat C be a category and \cat O a
(not necessarily small) full subcategory of \cat C. \Map C will
denote the category of maps of \cat C with commutative squares as
morphisms.
\begin{definition}
For any category \cat C, \emph{power category} $\Pow \cat C$ is the category
of all subsets of $\obj{\cat C}$ and \emph{enriched functions} as morphisms,
i.e., for any $A, B\subset \obj{\cat C}$ the set of morphisms from $A$ to $B$
is the set of all functions $F\colon A\rarrow \mor{\cat C}$, such that for any
$a\in A$, $\domain{F(a)} = a$ and $\range{F(a)}\in B$. Suppose $A, B, C\subset
\obj{\cat C}$. If $F_1\colon A\rarrow \mor{\cat C}$ is a morphism from $A$ to
$B$ and $F_2\colon B\rarrow \mor{\cat C}$ is a morphism from $B$ to $C$ in
$\Pow \cat C$, then their composition $F_3 = F_2\circ F_1$ is defined by
$F_3(a) = F_2(\range{F_1(a)}) \circ_{\cat C}F_1(a)$. For any set $A$, the
identity morphism $\id_A$ is a function that corresponds to every $a\in A$, the
identity morphism $\id_a$.
\end{definition}

\begin{definition}
Let $O$ be a class of objects in \cat C. For any object $c$ of \cat C, a set
of arrows $\psi_i\colon o_i\rarrow c$ has the \emph{factorization property}
with respect to $O$ if for every $o\in O$ and $\psi\colon o\rarrow c$ there is
a commutative triangle in \cat C
\[
\xymatrix@!0@ur{ o
\ar@{-->}[dr] \ar[d]_\psi\\
c       & o_i \ar[l]^{\psi_i}}
\]
for some $i$.
\end{definition}
The next definition has appeared in \cite[1.1]{Farjoun}.
\begin{definition}
A class $O$ of objects in a category \cat C is \emph{locally small} if for any
object $c\in \cat C$ there exists a set of arrows $\psi_i\colon o_i \rarrow c$
with $o_i\in O$, that has the factorization property with respect to $O$. A
full subcategory \cat O of a category \cat C is locally small if its class of
objects is locally small.
\end{definition}
\begin{example}
In the category of sets any class $S$ of pairwise non-isomorphic objects is locally small, since $S$ contains sets of cardinality bigger than any fixed cardinality $\kappa$, i.e. for
any set $X$ we can find a set $Y$ in $S$ which admits a surjective map $f\colon Y\rarrow X$. This map has the factorization property.
\end{example}
\begin{remark}
The condition for a class of maps to be locally small is dual to the classical
solution-set condition \cite{MacLane}.
\end{remark}

The previous definition admits a functorial version:

\begin{definition}\label{setup}
Let \cat O be a locally small subcategory of \cat C with a class of objects $O
= \obj{\cat O}$. A \emph{factorization setup} for the pair $(\cat C, O)$ is a
functor $F\colon\cat C\rarrow \Pow{(\Map C)}$ such that, for every
$c\in\obj{\cat C}$ and $f\in \mathrm{mor}_{\cat C}(c, c')$, we have
\begin{enumerate}
\item $F(c)$ is a set of maps of \cat C of the form $o\rarrow c$, $o\in O$;
\item $F(c)$ has the factorization property with respect to $O$;
\item $F(f)$ is a function that corresponds to any element $o \rarrow c$ of
$F(c)$ a commutative square (a morphism in $\Map C$) of the form
\[
\begin{CD}
 o  @>>>   o'\\
@VVV     @VVV\\
 c  @>f>>  c'
\end{CD}
\]
where $o' \rarrow c'$ is an element $F(c')$.
\end{enumerate}
\end{definition}

\begin{example}
In the category of groups the class of free groups has a factorization setup
$F$: for any group $G$, $F(G)$ is the set with one element
$\{\varepsilon_G\colon F_G \rarrow G\}$, where $\varepsilon_G$ is the
canonical map of the free group generated by the elements of $G$ onto $G$
(counit of the adjunction of the free functor with the forgetful functor); for
a homomorphism $\varphi\colon G_1 \rarrow G_2$ of groups, $F(\varphi)$ is
defined to be a function which assigns to $\varepsilon_{G_1}$ the following
commutative square:
\[
\begin{CD}
F_{G_1} @>{F_\varphi}>> F_{G_2}\\
 @VVV         @VVV\\
 G_1     @>{\varphi}>>  G_2
\end{CD}
\]
where $F_\varphi$ is the result of applying the free functor on the map of
underlying sets $\varphi\colon G_1 \rarrow G_2$. The factorization property of
$F(G)$ follows from the universal property of free groups.
\end{example}

We accept the following convention for the notation: if a class of objects in
a category \cat C is denoted by a capital letter, e.g.~$I$, then the full
subcategory of \cat C generated by the class $I$ is denoted by the script
letter \cat I and the factorization setup for the pair $(\cat C, I)$ is
denoted by the calligraphic letter $\cal I$.

\begin{example}\label{setup:set}
If a full subcategory \cat M of \cat C is small, i.e., $M=\obj{\cat M}$ is a
set, then there exists a factorization setup $\cal M\colon \cat C\rarrow
\Pow{(\Map C)}$, such that for every object $c$ in \cat C,
\[
\cal M(c) = \coprod_{m\in M} \mathrm{mor}_{\cat C}(m,c)
\]
and for any morphism $f\colon c\rarrow c'$ in \cat C, $\cal M(f)\colon \cal
M(c)\rarrow \cal M(c')$ is a function which corresponds to any element
$\varphi \colon m \rarrow c$ of $\cal M(c)$ the commutative square
\[
\begin{CD}
  m   @>{\id_m}>>  m\\
@V{\varphi}VV            @VV{f\varphi}V\\
  c   @>>f>      c'.
\end{CD}
\]
\end{example}

A simple property of factorization setups is given by the following
\begin{proposition}\label{disj:union}
Let $J$ and $K$ be two disjoint classes of objects in a category \cat C,
supplied with factorization setups \cal J and \cal K for the pairs $(\cat C,
J)$ and $(\cat C, K)$, respectively. Then the class of maps $L = J\sqcup K$ is
supplied with a factorization setup.
\end{proposition}
\begin{proof}
For any $c\in \cat C$ define $\cal L(c) = \cal J(c)\sqcup \cal K(c)$, and for
any morphism $f\colon c\rarrow c'$ in \cat C define $\cal L(f) = \cal
J(f)\sqcup \cal K(f)\subset \obj{\Map C}$. It is routine to check
Definition~\ref{setup}.
\end{proof}
\begin{corollary}\label{union}
If \cal J is a factorization setup for the pair $(\cat C, J)$ and $K$ is
a set of objects in \cat C, then the pair $(\cat C, J\cup K)$ may be
supplied with a factorization setup.
\end{corollary}
\begin{proof}
This follows from Proposition \ref{disj:union} applied for the disjoint
classes $J$ and $K'=K\setminus (K\cap J)$. A factorization setup for $K'$ is
provided by Example \ref{setup:set}.
\end{proof}

Let \cat C be a category, and let $I$ be a class of maps in \cat C. Recall
\cite{Hirschhorn, Hovey} that a map in \cat C is \emph{$I$-injective} if it
has the right lifting property with respect to every map in $I$; a map in \cat
C is \emph{$I$-projective} if it has the left lifting property with respect to
every map in $I$; a map in \cat C is \emph{$I$-cofibration} if it has the left
lifting property with respect to every $I$-injective map; a map in \cat C is
\emph{$I$-fibration} if it has the right lifting property with respect to
every $I$-projective map. The classes of $I$-injective maps, $I$-projective
maps, $I$-cofibrations and $I$-fibrations are denoted $I$-inj, $I$-proj,
$I$-cof and $I$-fib, respectively.

Suppose that \cat C contains all small colimits. A \emph{relative $I$-cell
complex} is a transfinite composition of pushouts of elements of $I$. We
denote the collection of relative $I$-cell complexes by $I$-cell.

Let \cat D be a collection of morphisms of \cat C, $A$ an object of \cat C and
$\kappa$ a cardinal. We say that $A$ is $\kappa$-\emph{small relative to \cat
D} if, for all $\kappa$-filtered ordinals $\lambda$ and all $\lambda$-sequences
\[
X_0\rarrow X_1\rarrow \cdots \rarrow X_\beta\rarrow \cdots,
\]
such that each map $X_\beta\rarrow X_{\beta+1}$ is in \cat D for $\beta+1 <
\lambda$, the natural map of sets
\[
\colim_{\beta<\lambda}\cat C(A, X_\beta)\rarrow \cat C(A,
\colim_{\beta<\lambda} X_\beta)
\]
is an isomorphism. We say that $A$ is \emph{small relative to} \cat D if it is
$\kappa$-small relative to \cat D for some $\kappa$. We say that $A$ is
\emph{small} if it is small relative to $\mor {\cat C}$.

Finally, we introduce the main notion of the current section: instrumented
class of maps.

\begin{definition}\label{instr}
Let \cat C be a category.  A locally small class $I$ of maps in \cat C is
called \emph{instrumented} (or equipped with an \emph{instrumentation}) if
\begin{enumerate}
\item there exists a cardinal $\kappa$ s.t. any $A\in\domain{I}$ is $\kappa$-small
        relative to $I$-cell,
\item there exists a factorization setup $\cal I$ for the pair $(\Map C, I)$ (see Definition~\ref{setup}).
\end{enumerate}
\end{definition}

\section{Example: Diagrams of spaces}\label{Diagrams}
Let $D$ be a small category. In this section we consider the equivariant model category on the $D$-shaped diagrams of spaces (see Definition~\ref{model-structure}) and show that the classes of generating cofibrations and generating trivial cofibrations (see below) may be equipped
with an instrumentation.

\subsection{Generating cofibrations and generating trivial cofibrations}
Let \cat O be the category of $D$-orbits. We define the class of generating
cofibrations to be $I = \{\dgrm T\otimes \partial \Delta^n \hookrightarrow
\dgrm T \otimes \Delta^n\}_{\dgrm T\in\cat O, n\geq 0}$ and the class of
generating trivial cofibrations to be $J = \{\dgrm T\otimes \Lambda^n_k
\tilde{\hookrightarrow} \dgrm T\otimes \Delta^n\}_{\dgrm T\in\cat O, n\geq k
\geq 0}$. Let \cat I and \cat J denote the full subcategories of $\Map{{\cal
S}^D}$ with classes of objects $I$ and $J$, respectively; $I$ and $J$ generate
the model structure of Definition~\ref{model-structure} in the sense that the
class of trivial fibrations equals $I$-inj and the class of fibrations equals
$J$-inj. The following proposition verifies for classes $I$ and $J$ the first
part of the definition of the instrumented class.

\subsection{Smallness of domains}
\begin{proposition}\label{small:domains}
Any $\dgrm T\otimes K \in \domain I \cup \domain J$, where \dgrm T is an orbit
and $K$ is a finite simplicial set, is $\aleph_0$-small with respect to both
the $I$-cell and $J$-cell.
\end{proposition}
\begin{proof}
Denote by $H$ the class of maps $\{\dgrm T\otimes K\rarrow \dgrm T\otimes L\}$, where $\dgrm T\in\cat O$ and $K\hookrightarrow L$ is a cofibration (inclusion) between finite simplicial sets. It will suffice to discuss colimits of $I$-cell maps along the first infinite ordinal $\omega$ only. Consider an $\omega$-sequence of $H$-cellular spaces $\emptyset=\dgrm
Z_0\rarrow\cdots\rarrow \dgrm Z_\alpha \rarrow\cdots \rarrow \dgrm Z_\omega$, where $\dgrm Z_\omega = \colim_{\alpha<\omega}Z_\alpha$.
Then by \cite[2.1]{Farjoun} the commutative square
\[
\begin{CD}
\dgrm Z_\alpha            @>i>>       \dgrm Z_\omega\\
    @VVV                                   @VVV\\
\colim_D \dgrm Z_\alpha   @>>>   \colim_D \dgrm Z_\omega
\end{CD}
\]
is a pullback for any finite ordinal $\alpha$.
\par
Let $\dgrm T\otimes K$ be the domain of some map in $H$ and $\dgrm T \otimes K
\overset k \rarrow \dgrm Z_\omega$ be any map. On the level of colimit the last
map factors through some finite stage $\alpha$, since $\colim_D \dgrm T\otimes
K = \ast\otimes K \cong K$ is $\omega$-small in $\cal S$, as a finite
simplicial set. Hence, by the universal property of a pullback, in the
commutative diagram of solid arrows
\[
\xymatrix{
\dgrm T\otimes K \ar[d]\ar@/^/[drr]^k \ar@{-->}[dr]_f \\
\colim_D \dgrm T\otimes K \cong K\ar@/_/[dr]      &         \dgrm Z_\alpha \ar[r]_i\ar[d]    &    \dgrm Z_\omega\ar[d]\\
               &         \colim_D \dgrm Z_\alpha \ar[r] &    \colim_D \dgrm Z_\omega
               }
\]
there exists the natural map $f\colon\dgrm T\otimes K\rarrow \dgrm Z_\alpha$,
such that $k=if$. We have obtained the required factorization of the original
map $k$.
\end{proof}

The following generalization of Proposition \ref{small:domains} is
required in Section \ref{loc}.
\begin{corollary}\label{smallness}
Any element of the set $\domain I \cup \domain J$ is $\aleph_0$-small with
respect to $I$-cof.
\end{corollary}
\begin{proof}
This follows from Proposition \ref{small:domains} and Proposition
\ref{Hirschhorn} below.
\end{proof}

\subsection{Factorization setup}
Before constructing a factorization setup for the locally small
classes of maps $I$ and $J$, a factorization setup will be
constructed on the locally small class of orbits.
\begin{proposition}
There exists a factorization setup for the pair $(\cal
S^D,\obj{\cat O})$.
\end{proposition}
\begin{proof}
We need to construct the functor $\cal O\colon \cal S^D\rarrow \Pow(\Map(\cal
S^D))$. Let \dgrm X be an object in $\cal S^D$, then define $\cal O(\dgrm X)$
to be a set of all maps of orbits into \dgrm X which are pullbacks of maps
$\ast \rarrow \colim_D \dgrm X$ along the canonical map $\dgrm X \rarrow
\colim_D \dgrm X$. The factorization property is satisfied by the universal
property of the pullback.
\par
If $\mor{\cal S^D}\ni f\colon\dgrm X\rarrow\dgrm Y$, then $\cal O(f)\colon \cal
O(\dgrm X) \rarrow \cal O(\dgrm Y)$ is a function which assigns to each
element $\{\dgrm T_x\rarrow \dgrm X \} \in \cal O(\dgrm X)$ the map in
$\Map(\cal S^D)$
\[
\xymatrix{
\dgrm T_x \ar[r]^F\ar[d]  &    \dgrm T_y\ar[d]\\
\dgrm X   \ar[r]^f        &    \dgrm Y, }
\]
where $x\in\colim\dgrm X$ is a point ($0$-simplex), $y=(\colim
f)(x)\in\colim\dgrm Y$, $\dgrm T_x$ and $\dgrm T_y$ are orbits over $x$ and
$y$, respectively, and the map $F$ is induced, naturally, by $f$.

It is routine to check that $\cal O$ is a functor.
\end{proof}
\begin{proposition}\label{factorization:setup}
There exist factorization setups $\cal I$ and $\cal J$ for the pairs
$(\Map(\cal S^D),\mhyphen I)$ and $(\Map(\cal S^D), J)$, respectively.
\end{proposition}
\begin{proof} We construct the factorization setup only for $I$; the
construction for $J$ is similar. Note that $\cal I\colon \Map(\cal S^D)
\rarrow \Pow{(\Map(\Map(\cal S^D)))}$.
\par
Let $f\colon\dgrm X\rarrow \dgrm Y$ be an object of $\Map(\cal S^D)$, then
define for each $n\geq 0$, $\dgrm W_{f,n}$ to be the pullback in the following
diagram:
\[
\xymatrix{
\dgrm W_{f,n} \ar[d] \ar[r]       &    \hom (\partial\Delta^n,\dgrm X) \ar[d] \\
\hom (\Delta^n,\dgrm Y) \ar[r]    &    \hom (\partial\Delta^n,\dgrm Y).
 }
\]
For $n$ fixed, let $I_{f,n}\subset \obj{\cat I\overcat f}$ which
consists of all commutative squares,
\[
\xymatrix{
\dgrm T\otimes\partial\Delta^n \ar[d] \ar[r]    &    \dgrm X \ar[d]^{f} \\
\dgrm T\otimes \Delta^n \ar[r]                  &    \dgrm Y }
\]
which correspond bijectively, by adjointness, to the commutative diagrams,
\[
\xymatrix{
\dgrm T  \ar@/_/[ddr] \ar@/^/[drr] \ar@{.>}[dr]|-{\varphi} \\
    & \dgrm W_{f,n} \ar[d] \ar[r]        & \hom(\partial\Delta^n,\dgrm X) \ar[d]^{f^{\partial\Delta^n}} \\
    & \hom(\Delta^n,\dgrm Y) \ar[r]  & \hom(\partial\Delta^n,\dgrm Y)
}
\]
where $\varphi$ runs through the \underline{set} $\cal O(\dgrm W_{f,n})$. Then
define $\cal I(f)= \bigcup_{n\geq 0}I_{f,n}$. The factorization property
holds by \cite[1.3, Lemma]{Farjoun}.
\par
For each map $\mor{\cal S^D}\ni g=(g_1,g_2)\colon f_1\rarrow f_2$, i.e., for
each commutative square
\[
\xymatrix{
\dgrm X_1 \ar[d]_{f_1} \ar[r]^{g_1}    &    \dgrm X_2 \ar[d]^{f_2} \\
\dgrm Y_1 \ar[r]_{g_2}                 &    \dgrm Y_2 }
\]
and for each $n\geq 0$, we obtain a natural map $\tilde{g}_n\colon \dgrm
W_{f_1,n}\rarrow \dgrm W_{f_2,n}$ between pullbacks in the following diagram:
\[
\xymatrix{
\dgrm W_{f_2,n} \ar[rrr] \ar[ddd] &&&  \hom (\partial\Delta^n,\dgrm X_2) \ar[ddd] \\
& \dgrm W_{f_1,n} \ar[d] \ar[r] \ar@{.>}[ul]_{\tilde g_n}^{\exists !}      &    \hom (\partial\Delta^n,\dgrm X_1) \ar[d] \ar[ur]\\
& \hom (\Delta^n,\dgrm Y_1) \ar[r] \ar[dl]   &    \hom (\partial\Delta^n,\dgrm Y_1) \ar[dr] \\
\hom (\Delta^n,\dgrm Y_2) \ar[rrr]  &&&  \hom
(\partial\Delta^n,\dgrm Y_2) }
\]
Then the morphism $\cal I(g)\colon\cal I(f_1)\rarrow \cal I(f_2)$ (in the
category $\Pow{(\Map(\Map(\cal S^D)))}$) is a function which assigns to each
element of $\cal I(f_1)$ the following commutative cube:
\[
\xymatrix@!{
    & \dgrm T_1\otimes\Delta^n \ar[rr]\ar'[d][dd]      & & \dgrm Y_1 \ar[dd] \\
\dgrm T_1\otimes\partial\Delta^n \ar[ur]\ar[rr]\ar[dd] & & \dgrm X_1 \ar[ur]\ar[dd]\\
    & \dgrm T_2\otimes\Delta^n  \ar'[r][rr]            & & \dgrm Y_2\\
\dgrm T_2\otimes\partial\Delta^n \ar[rr]\ar[ur]        & & \dgrm X_2 \ar[ur] }
\]
where the upper face is the respective element of $\cal I(f_1)$, the lower
face is an element of $\cal I(f_2)$, the right face is the initial map $g$, and
the vertical arrows of the left face are induced by the upper arrow in the
following commutative square:
\[ \xymatrix{
\dgrm T_1      \ar[r]\ar[d]            &    \dgrm T_2\ar[d]\\
\dgrm W_{f_1,n}  \ar[r]^{\tilde g_n}     &    \dgrm W_{f_2,n.}
 }
\]
This square equals $(\cal O(\tilde g_n))(\dgrm T_1\rarrow \dgrm W_{f_1,n})$.
The functoriality follows from the naturality of the construction.
\end{proof}
\begin{theorem}
The classes $I$ and $J$ are equipped with an instrumentation.
\end{theorem}
\begin{proof}
This follows from Proposition \ref{small:domains} and Proposition
\ref{factorization:setup}.
\end{proof}

The following result provides an initial motivation for the generalization of
Quillen's small object argument to the instrumented classes of maps, which we
postpone until the next section.
\begin{corollary}\label{functorial}
The factorizations in the model category of diagrams of spaces
generated by the full collection of orbits are functorial.
\end{corollary}
\begin{proof}
This follows from the generalized small object argument
\ref{argument} applied to the instrumented classes $I$ and $J$.
\end{proof}

\section{A generalization of Quillen's small object argument}\label{Quillen}
The main tool for the construction of the factorizations in the model
categories and localizations thereof is Quillen's small object argument
\cite{Hirschhorn, Hovey, Quillen}. However, in its original form, the argument
allows for the application in the cofibrantly generated categories only. We
propose here a generalization which may be used in a wider class of the model
categories.
\begin{proposition}[The generalized small object argument]\label{argument}
Suppose \cat C is a category containing all small colimits, and $I$ is an
instrumented class of maps in \cat C. Then there is a functorial factorization
$(\gamma,\delta)$ on \cat C such that, for all morphisms $f$ in \cat C, the
map $\gamma(f)$  is in $I$-cell and the map $\delta(f)$ is in $I$-inj.
\end{proposition}
\begin{proof}
Given the cardinal $\kappa$ such that every domain of $I$ is $\kappa$-small
relative to $I$-cell, we let $\lambda$ be a $\kappa$-filtered ordinal.
\par
To any map $f\colon X \rarrow Y$ we will associate a functor $Z^f\!\colon
\lambda \rarrow \cat C$ such that $Z_0^f = X$, and a natural transformation
$\rho^f\colon Z^f \rightarrow Y$ factoring $f$, i.e., for each $\beta <
\lambda$ the triangle
\[
\xymatrix@ur{
                                       X\ar[dr]^f \ar[d]\\
Z_f^\beta \ar[r]_{\rho^f_\beta}      & Y }
\]
is commutative. Each map $i^f_\beta\colon  Z_\beta^f \rightarrow
Z_{\beta+1}^f$ will be a pushout of a coproduct of maps of $\domain{\cal I(f)}
\subset I$, i.e., $i^f_\beta \in I$-cell.
\par
We will define $Z^f$ and $\rho^f\colon Z^f \rarrow Y$ by transfinite induction,
beginning with $Z^f_0 = X$ and $\rho_0^f = f$. If we have defined
$Z_{\alpha}^f$ and $\rho_{\alpha}^f$ for all $\alpha < \beta$ for some limit
ordinal $\beta$, define $Z_{\beta}^f = \colim_{\alpha<\beta} Z_{\alpha}^f$,
and  define $\rho_\beta^f$ to be the map induced, naturally, by the
$\rho_\alpha^f$. Having defined $Z_\beta^f$ and $\rho_\beta^f$, we define
$Z_{\beta+1}^f$ and $\rho_{\beta+1}^f$ as follows. The class $I$ of maps is
instrumented, hence $\cal I(\rho_\beta ^f)$ is the set of all commutative
squares of the following form:
\[
\begin{CD}
A  @>>> Z_\beta^f \\
@VgVV  @V{\rho_\beta^f}VV\\
B  @>>> Y
\end{CD}
\]
where $g \in \domain{\cal I(\rho_\beta^f)}\subset I$. For $s\in S=\cal I(
\rho_ \beta^f)$, let $g_s\colon A_s\rarrow B_s$ denote the corresponding map of
$I$. Define $Z_{\beta+1}^f$ to be the pushout in the diagram below,
\[
\begin{CD}
{\coprod_{s\in S}A_s}  @>>> Z_\beta^f \\
@V{\coprod g_s}VV               @VVV\\
{\coprod_{s\in S}B_s}  @>>> {Z_{\beta+1}^f.}
\end{CD}
\]
Define $\rho_{\beta+1}^f$ to be the map induced by $\rho_{\beta}^f$.
\par
For each morphism $g\colon f_1\rarrow f_2$ in the category \Map C, i.e., for
each commutative square
\[
\begin{CD}
{X_1}      @>{g^1}>>    X_2 \\
@V{f_1}VV            @V{f_2}VV\\
{Y_1}      @>{g^2}>>    {Y_2}
\end{CD}
\]
we define a natural transformation $\xi^g\colon Z^{f_1}\rarrow Z^{f_2}$ by
transfinite induction over small ordinals, beginning with $\xi^g_0 = g^1$. If
we have defined $\xi^g_\alpha$ for all $\alpha < \beta$ for some limit ordinal
$\beta$, define $\xi^g_\beta = \colim_{\alpha < \beta} \xi^g_\alpha$. Having
defined $\xi^g_\beta$, we define $\xi^g_{\beta+1}\colon Z^{f_1}_{\beta+1}
\rarrow Z^{f_2}_{\beta+1}$ to be the \emph{natural} map induced by $g_\beta =
(\xi^g_\beta, g^2)\colon  \rho^{f_1}_\beta \rarrow \rho^{f_2}_\beta$, namely
the \emph{unique} map between the pushouts of the horizontal lines of the
following diagram which preserves its commutativity:
\[
\begin{CD}
{\coprod_{s\in S}B_s}  @<<<  {\coprod_{s\in S}A_s}    @>>>    Z_\beta^{f_1} \\
      @VVV                       @VVV                           @VV{\xi_\beta^g}V\\
{\coprod_{t\in T}B_t}  @<<<  {\coprod_{t\in T}A_t}    @>>>    Z_\beta^{f_2}. \\
\end{CD}
\]
Here $S = \cal I(\rho^{f_1}_\beta)$ and $T = \cal I(\rho^{f_2}_\beta)$. The
first two vertical maps are induced, naturally, by the function $\cal
I(g_\beta)$. The commutativity of the above diagram follows readily from the
condition on the functor $\cal I$ to be a factorization setup.
\par
The required functorial factorization $(\gamma,\delta)$ is obtained when we
reach the limit ordinal $\lambda$ in the course of our induction. Then we
define $\gamma(f)\colon  X \rarrow Z^f_{\lambda}$ to be the (transfinite)
composition of the pushouts, and $\delta(f) = \rho^f_{\lambda}\colon
Z^f_{\lambda} \rarrow Y$. It follows from \cite[2.1.12, 2.1.13]{Hovey} that
$\gamma(f)$ is a relative $I$-cell complex.
\par
To complete the definition of the functorial factorization (see
\cite[1.1.1]{Hovey}, \cite[1.1.1]{Hovey_err}) we need to define for each
morphism $g\colon  f_1 \rarrow f_2$ a natural map $(\gamma,\delta)^g\colon
Z^{f_1}_{\lambda} \rarrow Z^{f_2}_{\lambda}$ which makes the appropriate
diagram commutative. Take $(\gamma,\delta)^g = \xi^g_\lambda$.
\par
It remains to show that $\delta(f) = \rho^f_{\lambda}$ has the right lifting
property with respect to $I$. To see this, suppose we have a commutative
square as follows:
\[
\begin{CD}
  C             @>h'>>              Z^f_\lambda \\
@VlVV                        @VV{\rho^f_\lambda}V\\
  D             @>k'>>              Y
\end{CD}
\]
where $l$ is a map of $I$. But the class $I$ is instrumented, hence there
exists a map $\domain{\cal I(\rho^f_\lambda)}\ni g\colon A\rarrow B$ such that
the map $(h',k')\in \mor{\Map C}$ factors through $g$. Hence it is enough to
construct a lift in the following diagram.
\[
\begin{CD}
  A          @>{h}>>        Z^f_\lambda \\
@VgVV                 @VV{\rho^f_\lambda}V\\
  B          @>{k}>>        Y
\end{CD}
\]
Due to the condition (1) of Definition \ref{instr} the object $A$ is
$\kappa$-small relative to $I$-cell, i.e., there is an ordinal $\beta <
\lambda$ such that $h$ is the composite $A\stackrel{h_\beta}{ \longrightarrow}
Z_\beta^f \longrightarrow Z^f_\lambda$. By construction, there is a map $B
\stackrel{k_\beta}{\longrightarrow} Z_{\beta+1}^f$ such that $k_\beta g =
i_\beta h_\beta$ and $k = \rho^f_{\beta+1}k_\beta$, where $i_\beta$ is the map
$Z_{\beta}^f\rarrow Z_{\beta+1}^f$. The composition $B \stackrel
{k_\beta}{\longrightarrow} Z_{\beta+1}^f\longrightarrow Z^f_\lambda$ is the
required lift in our diagram.
\end{proof}
The following results are straightforward generalizations of \cite[2.1.15,
2.1.16]{Hovey}, respectively. We record them for future reference.
\begin{corollary}\label{cof_retract}
Suppose $I$ is an instrumented class of maps in a category \cat C with all
small colimits. Then given $f\colon A \rarrow B$ in $I$-cof, there is a
$g\colon A \rarrow C$ in $I$-cell such that $f$ is a retract of $g$ by a map
that fixes $A$.
\end{corollary}
\begin{proposition}\label{Hirschhorn}
Suppose $I$ is an instrumented class of maps in a category \cat C that has all
small colimits. Suppose the domains of $I$ are small relative to $I$-cell, and
$A$ is some object that is small relative to $I$-cell. Then $A$ is in fact
small relative to $I$-cof.
\end{proposition}

\section{Properness}
In this section we show that the model category on the diagrams of spaces
generated by the collection of all orbits is proper. Recall that a model
structure is \emph{left proper} if weak equivalences are preserved under
pushouts along cofibrations. Dually, a model structure is \emph{right proper}
if weak equivalences are preserved under pullbacks along fibrations. A model
category is \emph{proper} if it is both left and right proper.

\subsection{Right properness} This is an immediate consequence of the right properness in the
category of simplicial sets.

In more details, given a pullback
\[
\xymatrix{
\dgrm W  \ar[r]\ar[d]     &    \dgrm X\ar[d]^\sim\\
\dgrm Z  \ar@{->>}[r]     &    \dgrm Y}
\]
and an orbit \dgrm T, the functor $\hom(\dgrm T, \cdot)$
preserves weak equivalences, fibrations and pullbacks. Therefore,
there is a pullback in the category of simplicial sets:
\[
\xymatrix{
\hom(\dgrm T, \dgrm W)  \ar[r]\ar[d]     &    \hom(\dgrm T, \dgrm X)\ar[d]^\sim\\
\hom(\dgrm T, \dgrm Z)  \ar@{->>}[r]     &    \hom(\dgrm T, \dgrm Y)}
\]
Hence, the right properness of $\cal S$ implies that the map $\hom(\dgrm
T, \dgrm W)\rarrow \hom(\dgrm T, \dgrm Z)$ is a weak equivalence of
spaces for any orbit \dgrm T, thus the map of diagrams $\dgrm W \rarrow
\dgrm Z$ is the weak equivalence of diagrams.

\subsection{Left properness}
The left properness is less straightforward. In fact, we provide different
proofs for the diagrams of simplicial sets and the diagrams of topological
spaces.

Let $\cal S$ be the category of simplicial sets.
\begin{proposition} \label{all:cofibrant}
Every diagram $\dgrm X\in \cal S^D$ is cofibrant in the equivariant model category.
\end{proposition}
\begin{proof}
Fix the diagram \dgrm X; we need to construct a lift in any commutative square of the form
\[
\xymatrix{
\emptyset  \ar[r]\ar[d]     &    \dgrm Y\ar@{->>}[d]^\sim\\
\dgrm X          \ar[r]     &    \dgrm Z}
\]
Let $\cat E$ be a small, full subcategory of the category of $D$-orbits, such
that \dgrm X is of orbit type $\cat E$ and all the free $D$-orbits are in
$\cat E$, i.e., there is the imbedding $D\hookrightarrow\cat E^{\op}$ that
associates to each $d\in \obj{D}$ the free orbit generated in $d$.
\par
Now we apply the functor $(\;\cdot\;)^{\cat E}$ on the commutative square above
and obtain a commutative square in the category of $\cat E^{\op}$-diagrams of
simplicial sets. Consider the Bousfield--Kan model category on the $\cat
E^{\op}$-diagrams; then the trivial fibration between $D$-diagrams becomes the
trivial fibration between $\cat E^{\op}$-diagrams. Since \dgrm X is of orbit
type $\cat E^{\op}$, $\dgrm X^{\cat E}$ is $\cat E^{\op}$-free by
\cite[3.7]{DF}. But in the Bousfield--Kan model category free diagrams of
simplicial sets are cofibrant, so there exists a lift in the square
\[
\xymatrix{
\emptyset  \ar[r]\ar@{^{(}->}[d]        &    \dgrm Y^{\cat E}\ar@{->>}[d]^\sim\\
\dgrm X^{\cat E} \ar[r] \ar@{-->}[ur]   &    {\dgrm Z}^{\cat E}
         }
\]
By Yoneda's lemma for any $d\in D$ and an arbitrary diagram \dgrm W, $\dgrm
W(d) \cong \dgrm W^{\cat E}(F^d)$. Hence, the lift in the above square
naturally `reduces' in the obvious sense to the lift in the original square.
\end{proof}
\begin{remark}
The last proposition settles affirmatively the conjecture stated in
\cite[2.3]{Farjoun}.
\end{remark}
\begin{corollary}\label{simplicial_left_proper}
The equivariant model category on $\cal S^D$ is left proper.
\end{corollary}
\begin{proof}
This follows from Proposition \ref{all:cofibrant} and
C.~L.~Reedy's theorem \cite{Reedy} (see \cite[13.1.2]{Hirschhorn}
for a modern exposition) which asserts that every pushout of a
weak equivalence between cofibrant objects along a cofibration is
a weak equivalence.
\end{proof}
\begin{proposition}\label{topological_left_proper}
The equivariant model category on the diagrams of topological spaces is left proper.
\end{proposition}
\begin{proof}
Suppose we are given a pushout in the category $\Top^D$ of diagrams of
topological spaces:
\[
\xymatrix{
\dgrm A  \ar[r]^f_\sim \ar@{^{(}->}[d]^g &    \dgrm B\ar[d]\\
\dgrm X  \ar[r]                          &    \dgrm Y }
\]
We have to show that the map $\dgrm X\rarrow\dgrm Y$ is a weak
equivalence.

By Corollary \ref{cof_retract} the cofibration $g\colon \dgrm A
\hookrightarrow \dgrm X$ is a retract of an $I$-cellular map $g'\colon \dgrm A
\hookrightarrow \dgrm X'$. Hence, the pushout of $f$ along $g$ is a retract of
the pushout of $f$ along $g'$ and it will suffice to show that the last map is
a weak equivalence.

The $I$-cellular cofibration $g'$ has a decomposition into a (transfinite)
sequence $\dgrm A = \dgrm A_0 \hookrightarrow \dgrm A_1 \hookrightarrow \cdots
\hookrightarrow \dgrm A_n \hookrightarrow \cdots \hookrightarrow\dgrm X'$, such
that $\dgrm A_{n+1}$ is a pushout of $\dgrm A_n$ along a map from $I$.
Therefore, the pushout of $f$ along $g$ is the colimit of the consecutive
pushouts of $f$ along the maps $g_n\colon  \dgrm A_n \hookrightarrow \dgrm
A_{n+1}$ for all $n$. Hence, by \mbox{Proposition \ref{small:domains}}, it is
enough to show that the pushout of a weak equivalence along $g_n$ is a weak
equivalence.

Any map in $I$ has a form $i\colon \dgrm T\otimes K \hookrightarrow\dgrm
T\otimes L$, where $K\hookrightarrow L$ is an inclusion of finite simplicial
sets and \dgrm T is an orbit. Consider the following commutative diagram in
which the left and right squares are pushouts:
\[
\xymatrix{
\dgrm T\otimes K \ar[r]
                 \ar@{^{(}->}[d]_i &
                 \dgrm A_n \ar[r]^{f_n}_\sim
                           \ar@{^{(}->}[d]_{g_n} &
                           \dgrm B_n\ar[d]\\
\dgrm T\otimes L \ar[r] &
                 \dgrm A_{n+1} \ar[r] &
                               \dgrm B_{n+1} }
\]
hence the outer square is a pushout.

By the axiom Q1, which was verified in Proposition \ref{Q1}, for any orbit
$\dgrm T'$, in the induced diagram of simplicial sets
\[
\xymatrix{
\hom(\dgrm T', \dgrm T\otimes K) \ar[rr]\ar@{^{(}->}[dd]_{\hom(\dgrm T',i)} &    & \hom(\dgrm T', \dgrm A_n) \ar[rr]^{\hom(\dgrm T',f_n)}_\sim \ar[dd]\ar@{^{(}->}[dl]&   &\hom(\dgrm T',\dgrm B_n)\ar[dd]\ar[dl]\\
                                                         & P_1 \ar@{-->}[dr] \ar@{-->}[rr]&                                                                     &P_2 \ar@{-->}[dr]&         \\
\hom(\dgrm T', \dgrm T\otimes L) \ar[rr] \ar[ur] \ar'[urr]!<-3 pt,-13 pt>[urrr]&    & \hom(\dgrm T', \dgrm A_{n+1}) \ar[rr]                               &   & \hom(\dgrm T', \dgrm B_{n+1})\\
 }
\]
the left and outer squares are, up to homotopy (i.e., up to a weak
equivalence), pushout diagrams. Consider the pushouts $P_1$ and $P_2$ of the
left and outer square respectively, i.e.,
\[
P_1 = \hom(\dgrm T', \dgrm T\otimes L) \coprod_{\hom(\dgrm T',
\dgrm T\otimes K)}\hom(\dgrm T', \dgrm A_n)
\]
and
\[
P_2 = \hom(\dgrm T', \dgrm T\otimes L) \coprod_ {\hom(\dgrm T',
\dgrm T\otimes K)}\hom(\dgrm T', \dgrm B_n),
\]
then the natural maps $P_1\rarrow \hom(\dgrm T', \dgrm A_{n+1})$ and $P_2
\rarrow \hom(\dgrm T', \dgrm B_{n+1})$ are weak equivalences of simplicial
sets.

The map $i$ in $\Top^D$ is a monomorphism, therefore the induced map
$\hom(\dgrm T',i)$ is a monomorphism (or a cofibration) of simplicial sets,
hence its cobase change $\hom(\dgrm T', \dgrm A_{n})\rarrow P_1$ is also a
cofibration.

Finally, \cite[7.2.14]{Hirschhorn} implies that $P_2= P_1 \coprod_ {\hom(\dgrm
T', \dgrm A_n)} \hom(\dgrm T', \dgrm B_n)$ and the map $P_1 \rarrow P_2$ is a
weak equivalence by the left properness of simplicial sets. Therefore, we can
conclude from the `2 out of 3' property that the map $\hom(\dgrm T', \dgrm
A_{n+1})\rarrow \hom(\dgrm T', \dgrm B_{n+1})$ is a weak equivalence of
simplicial sets and, hence, the original map $\dgrm A_{n+1} \rarrow \dgrm
B_{n+1}$ is a weak equivalence of diagrams of topological spaces.
\end{proof}

\section{Construction of the localization functor}\label{loc}
Let $S$ be a class of maps in $\cal S^D$. Without loss of generality, we may assume that the elements of $S$ satisfy the conditions listed in \ref{conditions} below. In this section we construct for such class $S$ a coaugmented functor $L_S\colon \cal S^D\rarrow \cal S^D$ such that for each $\dgrm X \in \cal S^D, \; L_S(\dgrm X)$ is $S$-local and the natural map $j_{\dgrm X}\colon \dgrm X\rarrow L_S(\dgrm X)$ is an $S$-equivalence (see below). We prove also that the natural map $j_{\dgrm X}\colon  \dgrm X\rarrow L_S\dgrm X$ is initial (up to homotopy) among all the maps of \dgrm X into $S$-local spaces, thus $L_S\dgrm X$ is characterized up to weak equivalence.
\par
Our construction is an extension of the classical constructions of localization functors with respect to a set of maps in a cofibrantly generated model category satisfying additional conditions \cite{Bous:factor}, \cite{Hirschhorn}. A useful summary of the classical construction is available at \cite[Appendix]{Javier}.

\subsection{Preliminaries on $S$-local spaces and $S$-equivalences}\label{conditions}
Throughout this section let us suppose that $S$ is a class of maps such that every $f\in S, \; f\colon \dgrm A \hookrightarrow \dgrm B$ is a cofibration between cofibrant diagrams. Assume, for simplicity, that all cofibrations in $S$ are non-trivial. We are going to construct the localization functor with respect to $S$ in this section, provided that $S$ satisfies the following conditions:
\begin{enumerate}
\item The class of horns of $S$ 
\[
\Hor(S)=\{\Delta^n\otimes \dgrm A \coprod_{\partial\Delta^n\otimes \dgrm A}
         \partial\Delta^n\otimes \dgrm B \rarrow \Delta^n\otimes \dgrm B \;|\; S \ni f\colon \dgrm A\rarrow \dgrm B, n\geq
         0\}
\]
may be equipped with an instrumentation with a factorization setup $\cal H$.
\item All the elements of $\domain{\Hor(S)}$ are $\kappa$-small with respect to cofibrations of  of diagrams for some fixed cardinal $\kappa$.
\end{enumerate}
These conditions are satisfied for example by any set of non-trivial cofibrations $S$. An example of a proper class of maps satisfying the conditions above will be given in Section~\ref{class:loc} below.

\begin{definition}\label{S-equivalence}
Let $S$ be a class of maps such that every $f\in S, \; f\colon \dgrm A \hookrightarrow \dgrm B$ is a cofibration between cofibrant diagrams.
\begin{itemize}
\item
A diagram \dgrm X is called $S$-\emph{local} if \dgrm X is fibrant and for
every $f\in S$, the induced map
\[
\hom(f,\dgrm X)\colon \hom(\dgrm B, \dgrm X) \rarrow \hom(\dgrm A, \dgrm X)
\]
is a weak equivalence of simplicial sets. If $S$ consists of the single map
$f\colon A\rarrow B$, then an $S$-local diagram will also be called $f$-local.
\item
A map $g\colon \dgrm C\rarrow \dgrm D$ is an $S$-\emph{local equivalence} (or
just an $S$-\emph{equivalence}) if for any cofibrant replacement $\tilde g$ of
$g$ and every $S$-local diagram \dgrm P the induced map
\[
\hom(\tilde{g},\dgrm P)\colon  \hom(\tilde{\dgrm D},\dgrm P) \rarrow
\hom(\tilde{\dgrm C},\dgrm P)
\]
is a weak equivalence of simplicial sets. If $S$ consists of the single map
$f\colon A\rarrow B$, then an $S$-local equivalence will also be called an
$f$-local equivalence (or an $f$-equivalence).
\end{itemize}
\end{definition}
\begin{remark}
Of course one needs to check that the notion of $S$-equivalence is
well-defined, i.e., it does not depend on the choice of the
cofibrant replacement. It follows from \cite[9.7.2]{Hirschhorn}.
We shall use also an $S$-local version of the Whitehead theorem
(see \cite[3.2.13]{Hirschhorn} for the proof).
\end{remark}

\begin{proposition}[$S$-local Whitehead theorem]\label{Whitehead}
A map $g\colon \dgrm Q_1 \rarrow \dgrm Q_2$ is a weak equivalence of $S$-local
spaces if and only if $g$ is an $S$-local equivalence.
\end{proposition}

\begin{proposition}\label{trivial}
A cofibration $g\colon \dgrm C \hookrightarrow \dgrm D$ of
cofibrant diagrams is an $S$-local equivalence if and only if $g$
has the homotopy left lifting property (see
\cite[9.4.2]{Hirschhorn} for the definition) with respect to all
the maps of the form $\dgrm Z\rarrow \ast$, where \dgrm Z is an
$S$-local diagram.
\end{proposition}
\begin{proof}
Follows immediately form the definitions and
\cite[9.3.1]{Hirschhorn}.
\end{proof}

\begin{proposition}\label{local}
A diagram ${\dgrm X}$ is $S$-local if and only if the map ${\dgrm X}\rarrow
\ast$ has the right lifting property with respect to the following families of
maps:
\begin{itemize}
\item generating trivial cofibrations $J$;
\item $\Hor(S)$.
\end{itemize}
\end{proposition}
\begin{proof}
The right lifting property for the members of $J$ is satisfied since $\dgrm X$
is fibrant by definition. For every element $f\in\Hor(S)$ it follows, by
adjunction, from the right lifting property of the map $\hom(f, \dgrm X)\colon
\hom(\dgrm B, {\dgrm X}) \twoheadrightarrow \hom(\dgrm A, {\dgrm X})$ with
respect to the generating cofibrations of simplicial sets
$\{\partial\Delta^n\hookrightarrow \Delta^n \;|\; n\geq 0\}$.
\end{proof}

\begin{proposition}\label{instr:K}
The class of maps $K = J\cup\Hor(S)$ may be equipped with an instrumentation.
\end{proposition}
\begin{proof}
The existence of a factorization setup \cal K for the pair $(\Map{{\calS} ^D},K)$ follows from Proposition~\ref{disj:union}, since we have chosen only non-trivial cofibrations in $S$, therefore the classes $J$ and $\Hor(S)$ are disjoint. The existence of the cardinal $\kappa$ such that the elements of $\domain{K}$ are $\kappa$-small relative to $K$-cell follows from the assumptions on $S$.
\end{proof}

\subsection{Construction of the functor $L_S$}
We construct the coaugmented functor $L_S$ by applying the generalized small
object argument, with respect to the instrumented (by Proposition
\ref{instr:K}) class of maps $K$, to factorize the map $\dgrm X \rarrow \ast$
into a $K$-cellular map, followed by a $K$-injective map. The obtained
functorial factorization
\[
\dgrm X \overset {j_{\dgrm X}} \longrightarrow L_S(\dgrm X)\longrightarrow \ast
\]
provides us with the coaugmented functor $L_S$, such that for any diagram
\dgrm X, the diagram $L_S(\dgrm X)$ is $S$-local by Proposition~\ref{local}.
\par
It remains to show that the natural coaugmentation map $j_{\dgrm X}$ is an
$S$-equivalence.
\begin{lemma}\label{lemma1}
Every map in $K$ is an $S$-equivalence.
\end{lemma}
\begin{proof}
Every map $g\in J$ is a trivial cofibration between cofibrant
diagrams, i.e., for any $S$-local (in particular fibrant) diagram
\dgrm Z, $\hom(g, \dgrm Z)$ is a trivial fibration of simplicial
sets. Hence, $g$ is an $S$-equivalence.
\par
It remains to show that every map in $\Hor(S)$ is an
$S$-equivalence. Every map in $\Hor(S)$ is a cofibration between
cofibrant objects, hence, by Proposition~\ref{trivial}, it is
enough to show that every map in $\Hor(S)$ has the \emph{homotopy}
left lifting property with respect to any map of the form $\dgrm Z
\rarrow \ast$, where \dgrm Z is $S$-local. This last property is
implied by \cite[9.4.8(1)]{Hirschhorn}.
\end{proof}
\begin{lemma}\label{lemma2}
A pushout of an $S$-equivalence $g$ that is also a cofibration between
cofibrant objects
\[
\xymatrix{
\dgrm A \ar@{^{(}->}[d]_g \ar [r]& \dgrm X \ar[d]\\
 \dgrm B \ar[r] & \dgrm Y
}
\]
in the equivariant model category on $\cal S^D$ is an $S$-equivalence again.
\end{lemma}
\begin{proof}
The following proof is a straightforward generalization of \cite[1.2.21]{Hirschhorn}. This argument significantly relies on the left properness.
\par
Factor the map $\dgrm A\rarrow\dgrm X$ as $\dgrm A \overset u \rarrow \dgrm
C\overset v \rarrow \dgrm X$, where $u$ is a cofibration and $v$ is a trivial
fibration. If we let $D$ be the pushout $\dgrm B\coprod_{\dgrm A}\dgrm C$,
then we have the commutative diagram
\[
\xymatrix{
\dgrm A \ar@{^{(}->}[d]_g \ar @{^{(}->}[r]^u & \dgrm C \ar@{^{(}->}[d]^k\ar@{->>}[r]^v_\sim & \dgrm X \ar[d]^h\\
 \dgrm B \ar[r]^s                            & \dgrm D \ar[r]^t                             & \dgrm Y
}
\]
in which $u$ and $s$ are cofibrations, and so \dgrm C and \dgrm D are
cofibrant. \cite[7.2.14]{Hirschhorn} implies that \dgrm Y is a pushout $\dgrm
X\coprod_{\dgrm C}\dgrm D$. Since $k$ is a cofibration and we are working in a
(left) proper model category (by Corollary \ref{simplicial_left_proper} and
Proposition \ref{topological_left_proper}), the map $t$ is a weak equivalence.
Thus, $k$ is a cofibrant approximation to $h$, and so it is sufficient to show
that $k$ induces a weak equivalence of mapping spaces to every $S$-local
diagram.

In any simplicial model category, a class of maps with the
homotopy left lifting property with respect to a map $p$ is closed
under pushouts \cite[9.4.9]{Hirschhorn}. Consider the collection
of maps $P = \{p_{\dgrm Z} \colon  \dgrm Z \rightarrow \ast \;|\;
\text{\dgrm Z is $S$-local} \}$. Then, by Proposition~
\ref{trivial}, a cofibration between cofibrant objects is an
$S$-local equivalence if and only if it has the homotopy left
lifting property with respect to any element of $P$. But the last
property is preserved under pushouts, hence $k$ in the diagram
above is an $S$-equivalence and so is $h$.
\end{proof}

\begin{lemma}\label{lemma3}
The class of $S$-equivalences between cofibrant diagrams in the simplicial
model category $\cal S^D$ is closed under coproducts (in the category of
maps), and the class of $S$-equivalences which are cofibrations is closed under
transfinite compositions in the left proper model category $\cal S^D$.
\end{lemma}
\begin{proof}
Let $g_\alpha$ be an $S$-equivalence for each $\alpha \in A$; then for any
$S$-local diagram \dgrm Z
\[
\hom(\coprod_{\alpha\in A}g_\alpha, \dgrm Z) = \prod_{\alpha\in
A}\hom(g_\alpha, \dgrm Z),
\]
where $\hom(g_\alpha, \dgrm Z)$ is a weak equivalence of simplicial sets,
therefore their product is also a weak equivalence.
\par
Let $\dgrm E_0\hookrightarrow\dgrm E_1\hookrightarrow\cdots
\hookrightarrow\dgrm E_\beta\hookrightarrow\cdots$ be a
$\lambda$-sequence of cofibrations which are also
$S$-equivalences. By \cite[17.9.4]{Hirschhorn} we may suppose,
without loss of generality, that all the diagrams $\dgrm E_i$ are
cofibrant. Then for each $S$-local diagram \dgrm Z there is a
$\lambda$-sequence of trivial fibrations of simplicial sets
\[
\hom(\dgrm E_0,\dgrm Z)\tilde {\twoheadleftarrow}\hom(\dgrm E_1,\dgrm Z)\tilde
{\twoheadleftarrow}\cdots \tilde {\twoheadleftarrow} \hom(\dgrm E_\beta,\dgrm
Z)\tilde {\twoheadleftarrow}\cdots.
\]
The inverse limit of the last sequence is a homotopy inverse limit, in
particular, the natural map $\hom(\dgrm E_0,\dgrm Z)\leftarrow \lim_{\beta<
\lambda}\hom(\dgrm E_\beta,\dgrm Z)= \hom(\colim_{\beta<\lambda} \dgrm
E_\beta,\dgrm Z)$ is a weak equivalence (compare to the inverse limit of the
constant tower). Thus, the natural map $\dgrm E_0\hookrightarrow \colim_{\beta
< \lambda} \dgrm E_\beta$ is an $S$-equivalence.
\end{proof}
\begin{proposition}
Every map in $K$-cell is an $S$-equivalence.
\end{proposition}
\begin{proof}
This follows from Lemmas \ref{lemma1}, \ref{lemma2} and \ref{lemma3}.
\end{proof}
\begin{corollary}\label{final}
The coaugmented functor $L_S$ is the $S$-localization functor.
\end{corollary}
\begin{proof}
By the construction of $L_S$, the natural map $j_{\dgrm X}\colon \dgrm X
\rarrow L_S{\dgrm X}$ is in $K$-cell, so, by the proposition above, $j_{\dgrm
X}$ is an $S$-equivalence. We postpone the proof of universality of $L_S$ until
the next section.
\end{proof}
\begin{lemma}\label{lemma4}
For any diagram \dgrm X, either cofibrant or not, the coaugmentation morphism
$j_{\dgrm X}\colon \dgrm X\rarrow L_S\dgrm X$ satisfies: for any $S$-local
diagram \dgrm P,
\[
\hom(j_{\dgrm X},\dgrm P) \colon  \hom(L_S\dgrm X, \dgrm P)\overset \simeq
\rarrow \hom(\dgrm X, \dgrm P)
\]
is a weak equivalence of simplicial sets.
\end{lemma}
\begin{proof}
Corollary~\ref{final} implies that $j_{\dgrm X}$ is an
$S$-equivalence, hence the result follows from
Definition~\ref{S-equivalence} and \cite[13.2.2(1)]{Hirschhorn}.
\end{proof}
\subsection{Universality and other properties of $L_S$}
\begin{proposition}[$L_S$ is initial]\label{initial}
For any map $g\colon \dgrm X \rarrow \dgrm P$ into an $S$-local diagram there
exists a factorization $\dgrm X \rarrow L_S\dgrm X\rarrow \dgrm P$ which is
unique up to simplicial homotopy.
\end{proposition}
\begin{proof}
By the construction of $L_S$, the natural map $j_{\dgrm X}\colon \dgrm X
\rarrow L_S{\dgrm X}$ is in $K$-cell. By Proposition \ref{local} the map $\dgrm
P \rarrow \ast$ is in $K$-inj. Then in the diagram
\[
\xymatrix{
\dgrm X  \ar[r]^g \ar@{^{(}->}[d]_{K\text{-cell}\;\ni j_{\dgrm X}} & \dgrm P\ar@{->>}[d]^{\in K\text{-inj}}\\
L_S\dgrm X \ar[r] \ar@{-->}[ur]                                    & \ast }
\]
the lift exists and provides the required factorization.

To show the uniqueness up to homotopy of the factorization above,
consider the map $\hom(j_{\dgrm X},\dgrm P) \colon  \hom(L_S\dgrm
X, \dgrm P)\overset \simeq \longrightarrow \hom(\dgrm X, \dgrm
P)$, which is a weak equivalence by Lemma \ref{lemma4}, since
\dgrm P is $S$-local. Then, by \cite[9.5.10]{Hirschhorn},
$j_{\dgrm X}$ induces a bijection of the sets of simplicial
homotopy classes of maps $j_{\dgrm X}^\ast\colon [L_S\dgrm X,
\dgrm P]\cong [\dgrm X, \dgrm P]$. The uniqueness, up to
simplicial homotopy, of the lifting follows from the injectivity
of the map $j_{\dgrm X}^\ast$.
\end{proof}
\begin{remark}
By \cite[II.2.5]{Quillen}, if two maps are simplicially homotopic, then they
are both left and right homotopic.
\end{remark}
\begin{corollary}
The natural maps $j_{L_S\dgrm X}, L_S(j_{\dgrm X})\colon  L_S\dgrm X
\rightrightarrows L_SL_S\dgrm X$ are simplicially homotopic weak equivalences.
\end{corollary}
\begin{proof}
The naturality of $j$ implies that the square
\[
\begin{CD}
\dgrm X            @>{j_{\dgrm X}}>>        L_S\dgrm X\\
@V{j_{\dgrm X}}VV                             @VV{j_{L_S\dgrm X}}V\\
L_S\dgrm X         @>>{L_S(j_{\dgrm X})}>   L_SL_S\dgrm X
\end{CD}
\]
is strictly commutative. Two possible paths provide two different
factorizations of the map $\dgrm X \rarrow L_SL_S\dgrm X$. By Proposition
\ref{initial} the natural maps $j_{L_S\dgrm X}, L_S(j_{\dgrm X})$ are
simplicially homotopic. Proposition \ref{Whitehead} implies that the map
$j_{L_S\dgrm X}$ is a weak equivalence; so does $L_S(j_{\dgrm X})$.
\end{proof}
\begin{proposition}[Divisibility]\label{divisibility}
For any two maps $g,h\colon L_S\dgrm X\rightrightarrows \dgrm P$ into any
$S$-local diagram \dgrm P, one has $h \overset s \sim g$ if and only if
$h\circ j_{\dgrm X} \overset s \sim g\circ j_{\dgrm X}$.
\end{proposition}
\begin{proof}
Suppose $h \overset s \sim g$. The simplicial homotopy between the maps $h$ and
$g$ is a $1$-simplex in $\hom(L_S\dgrm X, \dgrm P)$, so its image under
$j_{\dgrm X}^\ast$ provides a simplicial homotopy between $h\circ j_{\dgrm X}$
and $g\circ j_{\dgrm X}$.

Conversely, if $h\circ j_{\dgrm X} \overset s \sim g\circ j_{\dgrm X}$, then
the injectivity of the map $j_{\dgrm X}^\ast\colon [L_S\dgrm X, \dgrm P] \cong
[\dgrm X, \dgrm P]$ implies that $h \overset s \sim g$.
\end{proof}

\begin{proposition}[No zero divisors] \label{no-zero-divisors}
Suppose \dgrm W is a retract of $L_S\dgrm X$ for some \dgrm X. If the
composition $X\rarrow L_S\dgrm X\rarrow \dgrm W$ is null homotopic (see
Appendix~ \ref{null-homotopy} for the definition), then $\dgrm W\simeq \ast$.
\end{proposition}
\begin{proof}
First notice that the diagram \dgrm W is $S$-local (in particular, fibrant) as
a retract of the $S$-local diagram $L_S\dgrm X$. Suppose that the composition
$X\rarrow L_S\dgrm X\rarrow \dgrm W$ is null homotopic; then there exists, by
Proposition~ \ref{null}, a fibrant contractible diagram \dgrm U such that the
following solid arrow diagram commutes:
\[
\xymatrix{
\dgrm X  \ar[dr] \ar[r]^{j_{\dgrm X}}  & L_S\dgrm X\ar[r]_{r}\ar@{-->}[d]_p     & \dgrm W\ar@{.>}@/_1pc/[l]|i\\
                                       & \dgrm U\ar[ur]_q &\ar @{} [l]+<14pt,0pt>*{\simeq \dgrm\ast.} }
\]
The dashed arrow $p$ exists by the universal property of
Proposition~\ref{initial} and makes the left triangle commutative. By the
divisibility property of Proposition \ref{divisibility}, $r \overset s \sim
q\circ p$, since, upon precomposing with $j_{\dgrm X}$, these two maps are
equal.

Recall that \dgrm W is a retract of $L_S\dgrm X$, hence $\id_{\dgrm W} =
r\circ i \overset s \sim q\circ (p\circ i)$. Therefore, \dgrm W is a
\emph{homotopy} retract of \dgrm U.

Now apply the general machinery for detection of weak equivalences
\cite{Hirschhorn}: the map $q\colon \dgrm U\rarrow \dgrm W$ between two
fibrant spaces is a weak equivalence if and only if for any cofibrant \dgrm A
the induced map on simplicial homotopy classes $q_\ast\colon [\dgrm A,\dgrm
U]\rarrow [\dgrm A,\dgrm W]$ is a (natural) bijection. But $[\dgrm A,\dgrm U]
= \ast$, as \dgrm U is fibrant and contractible, and $[\dgrm A,\dgrm W]$ is a
retract of $[\dgrm A,\dgrm U]$, i.e., $[\dgrm A,\dgrm W] = \ast$. Hence, $\dgrm
W \simeq \dgrm U \simeq \ast$.
\end{proof}
\begin{corollary}
$L_S(\dgrm U)\simeq\dgrm\ast$ for any contractible diagram \dgrm U and $L
_S(\textup{nullmap})$ is a null map.
\end{corollary}
\begin{proof}
$L_S(\dgrm U)$ is a retract of $L_S(\dgrm U)$ by the identity morphism.
Moreover, $\id_{L_S(\dgrm U)} \circ j_{\dgrm U}\colon  \dgrm U\rarrow L_S(\dgrm
U)$ is null, since \dgrm U itself is contractible. Hence, by the proposition
above $L_S(\dgrm U) \simeq\dgrm\ast$.
\par
The second property follows immediately from the first one.
\end{proof}

\begin{proposition}\label{S-equiv}
A map $g\colon  \dgrm X \rarrow \dgrm Y$ is an $S$-local equivalence if and
only if $L_S(g)\colon  L_S \dgrm {X}\rarrow L_S \dgrm {Y}$ is a weak
equivalence.
\end{proposition}
\begin{proof}
In the commutative diagram
\[
\begin{CD}
    \dgrm {X}             @>{g}>>           \dgrm {Y}\\
 @V{j_ {\dgrm X}}VV               @VV{j_{\dgrm Y}}V\\
    L_S\dgrm {X}      @>>{L_S(g)}>       L_S\dgrm {Y}
\end{CD}
\]
the vertical arrows are $S$-local equivalences by construction.
Hence, the map $g$ is an $S$-local equivalence if and only if
$L_S(g)$ is an $S$-local equivalence by the `2 out of 3' property
of $S$-equivalences \cite[3.2.3]{Hirschhorn}. But $L_S(g)$ is a
map between two $S$-local spaces, hence $L_S(g)$ is an $S$-local
equivalence if and only if $L_S(g)$ is a weak equivalence by
Proposition \ref{Whitehead}.
\end{proof}

The coaugmentation map $j_{\dgrm X}$ is a cofibration for any diagram \dgrm X,
hence the subcategory of cofibrant diagrams is stable under localizations.
Consider the restriction of $L_S$ to the subcategory of cofibrant objects (do
nothing for the diagrams of simplicial sets) and denote the new functor by
$L_S^r$. Then $L_S^r$ is terminal with respect to $S$-local equivalences. In
more detail, we have the following
\begin{proposition}[$L^r_S$ is terminal]
On the subcategory of cofibrant objects the coaugmentation map $j_{\dgrm
X}\colon \dgrm X\rarrow L_S\dgrm X = L_S^r\dgrm X$ is terminal, up to homotopy,
among all $S$-local equivalences, i.e., for any $S$-equivalence of cofibrant
diagrams  $g\colon \dgrm X\rarrow \dgrm Y$, there exists an extension $\dgrm
X\rarrow \dgrm Y\overset l \rarrow L_S^r\dgrm X$ that is unique up to
(simplicial) homotopy with $l\circ g\overset s \sim {j_{\dgrm X}}\colon \dgrm
X\rarrow L_S^r\dgrm X$.
\end{proposition}
\begin{proof}
By Proposition \ref{S-equiv}, $L_S^r(g) = L_S(g)$ is a weak equivalence.
Moreover, this is a weak equivalence between two objects which are both
fibrant and cofibrant, so in the following commutative diagram of solid arrows
\[
\xymatrix{
   \dgrm X\ar[r]^g \ar[d]_{j_{\dgrm X}}      &    \dgrm Y\ar[d]_{j_{\dgrm Y}}\\
L_S\dgrm X\ar[r]^\sim_{L_S g}                & L_S\dgrm Y \ar@/^15pt/@{-->}[l]^q\\
}
\]
the map $L_S(g)$ has a simplicial homotopy inverse $q$ (all
notions of homotopy of maps between objects which are fibrant and
cofibrant coincide). Define $l=qj_{\dgrm Y}$; then, using the
commutativity of the diagram above and \cite[9.5.4]{Hirschhorn},
we obtain $lg = qj_{\dgrm Y}g = (qL_S(g))j_{\dgrm X} \overset s
\sim \id_{L_S(\dgrm X)}j_{\dgrm X} = j_{\dgrm X}$.

To show the uniqueness up to simplicial homotopy of $l$, suppose
there exists another map $l'\colon  \dgrm Y \rarrow L_S\dgrm X$
such that $l'g \overset s \sim j_{\dgrm X}$. Then $l'$ factors
through $L_S\dgrm Y$ (since $L_S$ is initial), i.e. there exists
an arrow $q'\colon L_S\dgrm Y\rarrow L_S\dgrm X$ such that $l'= q'
j_{\dgrm Y}$. It will suffice to show that $q' \overset s \sim q$
since \cite[9.5.4]{Hirschhorn} implies that $l' \overset s \sim
l$. By assumption, $l'g \overset s \sim j_{\dgrm X}$, so
$q'j_{\dgrm Y}g \overset s \sim j_{\dgrm X}$ or, equivalently,
$q'L_S(g) j_{\dgrm X} \overset s \sim \id_{L_S{ \dgrm X}}j_{\dgrm
X}$. By divisibility, $q'L_S(g) \overset s \sim \id_{L_S{ \dgrm
X}}$. But $L_S(g)$ is a weak equivalence of cofibrant objects,
therefore, by \cite[9.5.12]{Hirschhorn}, $L_S(g)^\ast\colon
[L_S\dgrm Y, L_S\dgrm X] \overset \approx \rarrow [L_S\dgrm X, L_S
\dgrm X]$ is a bijective map of simplicial homotopy classes which
satisfies $L_S(g)^\ast([q]) = L_S(g)^\ast ([q'])$, hence
$[q]=[q']$ or $q \overset s \sim q'$.
\end{proof}

\section{Fixed-pointwise localization --- Localization with respect to a class of maps}\label{class:loc}
In the previous section we developed the localization theory of diagrams of spaces with respect to a class of maps of diagrams subject to certain conditions \ref{conditions}. But a large part of the equivariant homotopy theory \cite{Bredon, DF, DZ} uses the `fixed-pointwise' approach and its generalizations, so it is natural to ask whether for any map $f\colon A\hookrightarrow B$ of simplicial sets there exists a localization functor $L$, which induces $f$-equivalences of fixed-point sets $\hom(\dgrm T, \dgrm X)\rarrow \hom(\dgrm T, L\dgrm X)$ for each orbit \dgrm T. We do not know in general whether it is possible to find a set of maps of diagrams such that the localization with respect to it gives the required functor $L$. But in one simple case $f\colon \emptyset \cofib \ast$, discussed in the companion paper \cite{PhDII}, we know that this is impossible. In this section we show how to apply the generalized small object argument to the localizations with respect to certain classes of maps.

There is no point in considering fixed-pointwise localizations with respect to a set of maps since, in the category of spaces, the localization with respect to any set of maps is equivalent to the localization with respect to a single map: take this single map to be the coproduct of the set of maps, in case that there is no map of the form $\emptyset\hookrightarrow X$, $X\neq \emptyset$; otherwise, this map $\emptyset\hookrightarrow X$ will induce the same localization functor as the whole set.

Fixed-pointwise localization with respect to homology in the category of spaces with a compact Lie group action was constructed in \cite{MMT}. Our construction is new only for the diagram shapes which lead to non-cofibrantly generated model structures on $\cal S^D$ \cite{Chorny}, otherwise it is covered by the classical localization framework.

Given a non-trivial cofibration $f\colon A\hookrightarrow B$ of simplicial sets, consider the following class of maps $F = \{f\otimes \dgrm T\colon A \otimes \dgrm T \hookrightarrow B\otimes \dgrm T \;|\; \dgrm T\in \cat O\}$. Then a diagram \dgrm Z is $F$-local if and only if it is fibrant and for each orbit \dgrm T the space of `\dgrm T-fixed points', $\hom(\dgrm T, \dgrm Z)$, is $f$-local:
\begin{itemize}
\item
\dgrm Z is fibrant $\Leftrightarrow \hom(\dgrm T, \dgrm Z)$ is fibrant for each
$\dgrm T \in \cat O$, by definition;
\item
$\hom(f\otimes \dgrm T, \dgrm Z)\colon \hom(B \otimes \dgrm T, \dgrm Z) \rarrow
\hom( A \otimes \dgrm T, \dgrm Z)$ is a weak equivalence $\Leftrightarrow
\hom(B, \hom (\dgrm T, \dgrm Z)) \rarrow  \hom(A, \hom (\dgrm T, \dgrm Z))$ is
a weak equivalence for each $\dgrm T\in \cat O$, by adjunction.
\end{itemize}

\begin{example}\label{first-example}
Let $f\colon \emptyset \hookrightarrow \ast$ be a map in $\cal S$; then $F =
\{ f\otimes \dgrm T \colon \emptyset \hookrightarrow \dgrm T \;|\; \dgrm T\in
\cat O \}$. The $f$-localization functor on the category of spaces assigns to
any space $X$ a contractible space $L_f(X)$, for a space is $f$-local iff it
is fibrant and contractible. By the considerations above, a diagram \dgrm X is
$F$-local iff its \dgrm T-fixed-point space is $f$-local for any orbit \dgrm
T, i.e., fibrant and contractible, hence \dgrm X is $F$-local iff it is
fibrant and contractible.

Any map of diagrams is an $F$-equivalence, likewise any map of spaces is an
$f$-equivalence. In other words, a map of diagrams $g\colon \dgrm X\rarrow
\dgrm Y$ is an $F$-local equivalence if and only if the induced map of
fixed-point spaces $\hom(\dgrm T, g) \colon \hom(\dgrm T, \dgrm X) \rarrow
\hom(\dgrm T, \dgrm Y)$ is an $f$-equivalence for each orbit $\dgrm T\in \cat
O$.
\end{example}

The example above has the following generalization:

\begin{proposition}\label{F-local-equiv}
A map $g\colon \dgrm X\rarrow \dgrm Y$ is an $F$-local equivalence if and only
if for each orbit $\dgrm T\in \cat O$ the map $\hom(\dgrm T, g) \colon
\hom(\dgrm T, \dgrm X) \rarrow \hom(\dgrm T, \dgrm Y)$ is an $f$-equivalence
of simplicial sets.
\end{proposition}
\begin{proof}
Let $\cat E \subset \cat O$ be a full small subcategory of the category of
orbits such that \dgrm X and \dgrm Y are of orbit type \cat E. The set of
orbits $\obj{\cat E}$ is a set of orbits in the sense of \cite{DK} (see
\cite{Farjoun} for the proof), and induces a simplicial model structure on the
category of $D$-diagrams which is Quillen equivalent to the Bousfield--Kan
model category on ${\cat E}^{\op}$-shaped diagrams of spaces \cite{DK}.
Moreover, these model categories are cofibrantly generated, therefore they have
functorial factorizations, and it was shown in \cite{DK2} that their
simplicial homotopy categories are homotopy equivalent.

A map $g\colon  \dgrm X\rarrow \dgrm Y$ is an $F$-local equivalence if for any
cofibrant replacement (in the model category generated by \underline{all}
orbits) $\tilde g\colon  \tilde{\dgrm X}\rarrow \tilde{\dgrm Y}$ and any
$F$-local diagram \dgrm Z the induced map on function complexes $\hom( \tilde
g, \dgrm Z)\colon  \hom(\tilde{\dgrm Y}, \dgrm Z)\rarrow \hom(\tilde{\dgrm X},
\dgrm Z)$ is a weak equivalence of simplicial sets. By the construction (in
Section~ \ref{Diagrams}) of the cofibrant replacement (which is a
$D$-$CW$-complex of the same orbit type as the original space), $\tilde{\dgrm
X}$ and $\tilde {\dgrm Y}$ are cofibrant in the model category generated by
the objects of \cat E and \dgrm Z is fibrant in both model categories, hence
the function complexes $\hom(\tilde{\dgrm X}, \dgrm Z)$ and $\hom(\tilde{\dgrm
Y}, \dgrm Z)$ are weakly equivalent to the homotopy function complexes of the
simplicial homotopy category \cite{DK2}. But $(\dgrm X)^{\cat E}$ and $(\dgrm
Y)^{\cat E}$ are cofibrant in the Bousfield--Kan model category as free
diagrams over ${\cat E}^{\op}$, and $(\dgrm Z)^{\cat E}$ is fibrant, hence the
function complexes $\hom((\tilde{\dgrm X})^{\cat E}, (\dgrm Z)^{\cat E})$ and
$\hom((\tilde{\dgrm Y})^{\cat E}, (\dgrm Z)^{\cat E})$ are also weakly
equivalent to the homotopy function complexes of the corresponding simplicial
homotopy category.

The homotopy equivalence between simplicial homotopy categories implies that
the map $\hom(\tilde g, \dgrm Z)\colon  \hom(\tilde{\dgrm Y}, \dgrm Z)\rarrow
\hom (\tilde{\dgrm X}, \dgrm Z)$ is a weak equivalence if and only if the map
$\hom( (\tilde g)^{\cat E}, (\dgrm Z)^{\cat E})\colon  \hom((\tilde{\dgrm Y})^
{\cat E}, (\dgrm Z)^{\cat E})\rarrow \hom((\tilde{\dgrm X})^{\cat E}, (\dgrm
Z)^{\cat E})$ is a weak equivalence. But $\tilde{\dgrm X}$ and $\tilde{\dgrm
Y}$ are weakly equivalent to \dgrm X and \dgrm Y, hence $(\tilde{\dgrm
X})^{\cat E} \cong ({\dgrm X})^{\cat E}$ and $(\tilde{\dgrm Y})^{\cat E} \cong
({\dgrm Y})^{\cat E}$ are weak equivalences of cofibrant (free) objects in the
Bousfield--Kan model category. Therefore, the maps $\hom(\tilde g, \dgrm Z)$
and $\hom((\tilde g)^{\cat E}, (\dgrm Z)^{\cat E})$ are weak equivalences if
and only if $\hom((g)^{\cat E}, (\dgrm Z)^{\cat E})\colon  \hom(({\dgrm Y})^
{\cat E}, (\dgrm Z)^{\cat E})\rarrow \hom(({\dgrm X})^{\cat E}, (\dgrm
Z)^{\cat E})$ is a weak equivalence. Compare \cite[5.13]{DF}.

Now we can prove the proposition. Suppose that all maps induced by the map $g$
on fixed-point sets are $f$-equivalences. Let \dgrm Z be any $F$-local diagram.
The equivariant function complex $\hom(({\dgrm X})^{\cat E}, (\dgrm Z)^{\cat
E})$ may be represented as a homotopy inverse limit over the \emph{twisted
arrow category} $a{\cat E}^{\op}$ (objects of $a{\cat E}^{\op}$ are morphisms
of ${\cat E}^{\op}$ and arrows of $a{\cat E}^{\op}$ are commutative squares
\[
\begin{CD}
 e_0 @<<< e_0'\\
@VVV      @VVV\\
 e_1 @>>> e_1'
\end{CD}
\]
in ${\cat E}^{\op}$) \cite[3.3]{DK1}. Since $(\dgrm X)^{\cat E}$ is cofibrant
(because it is free) and $(\dgrm Z)^{\cat E}$ is fibrant, the following map is
a weak equivalence:
\[
\hom(({\dgrm X})^{\cat E}, (\dgrm Z)^{\cat E}) \overset \sim \longrightarrow
\holim_{a\cat E^{\op}} \hom_a(({\dgrm X})^{\cat E}, (\dgrm Z)^{\cat E}),
\]
where $\hom_a(({\dgrm X})^{\cat E}, (\dgrm Z)^{\cat E})$ is an $a\cat
E^{\op}$-diagram of simplicial sets in which for each $(e_0\rarrow e_1)\in
a\cat E^{\op}$ there is assigned the simplicial set $\hom(({\dgrm X})^{\cat
E}(e_0), (\dgrm Z)^{\cat E}(e_1))$.

By assumption, the induced map between $a\cat E^{\op}$-diagrams
\[
\hom_a(g^{\cat E}, (\dgrm Z)^{\cat E})\colon  \hom_a((\dgrm Y)^{\cat E}, (\dgrm
Z)^{\cat E})\rarrow \hom_a((\dgrm X)^{\cat E}, (\dgrm Z)^{\cat E})
\]
is an objectwise weak equivalence, since each entry of the diagram $(\dgrm
Z)^{\cat E}$ is an $f$-local space. Hence, the induced map on the homotopy
inverse limits is a weak equivalence $\hom((g)^{\cat E}, (\dgrm Z)^{\cat
E})\colon  \hom(({\dgrm Y})^ {\cat E}, (\dgrm Z)^{\cat E})\overset \sim
\longrightarrow \hom(({\dgrm X})^{\cat E}, (\dgrm Z)^{\cat E})$. Then, by the
discussion above, the map $\hom(\tilde g, \dgrm Z)\colon  \hom(\tilde{\dgrm Y},
\dgrm Z)\rarrow \hom (\tilde{\dgrm X}, \dgrm Z)$ is a weak equivalence for any
$F$-local diagram \dgrm Z, i.e., the map $g$ is an $F$-local equivalence.

Alternatively, one can think of an equivariant function complex as an end. Then
this is a homotopy end, since $(\dgrm X)^{\cat E}$ is free, and the same
conclusion follows from the results in \cite{DDF}.

Suppose now that $g$ is an $F$-local equivalence. We need to show that for each
orbit \dgrm T the induced map on the \dgrm T-fixed-point space is an
$f$-equivalence. Let $V_{\dgrm T}(\dgrm T')=\hom_{\cat E^{\op}}(\dgrm T',
\dgrm T)$ be a diagram of simplicial sets over $\cat E$; then by the dual of
Yoneda's lemma we obtain the natural isomorphism $(\dgrm X)^{\cat E} \mathbin
{\otimes_ {\cat E^{\op}}} V_{\dgrm T} \cong (\dgrm X)^{\cat E}(\dgrm T) \cong
\hom(\dgrm T, \dgrm X)$. Take $W$ to be any $f$-local simplicial set. Then in
the commutative diagram
\[
\xymatrix{
   \hom((\dgrm Y)^{\cat E}(\dgrm T), W) \ar[r]^\cong \ar[d]_{\hom(g^{\cat E}|_{\dgrm T},}^{W)}    &    \hom((\dgrm Y)^{\cat E}\mathbin {\otimes_ {\cat E^{\op}}} V_{\dgrm T}, W)  \ar[r]^\cong \ar[d]     &  \hom((\dgrm Y)^{\cat E}, \hom(V_{\dgrm T}, W)) \ar[d]_{\hom(g^{\cat E},}^{\hom(V_{\dgrm T}, W))} \\
   \hom((\dgrm X)^{\cat E}(\dgrm T), W) \ar[r]^\cong                                              &    \hom((\dgrm X)^{\cat E}\mathbin {\otimes_ {\cat E^{\op}}} V_{\dgrm T}, W)  \ar[r]^\cong            &  \hom((\dgrm X)^{\cat E}, \hom(V_{\dgrm T}, W)), \\
}
\]
where $\hom(V_{\dgrm T}, W))$ is an $\cat E^{\op}$-diagram of $f$-local spaces
\cite[A.8(e.2)]{Farjoun-book}, the left vertical arrow is a weak equivalence
if and only if the right vertical arrow is a weak equivalence. But the diagram
$\hom(V_{\dgrm T}, W)$ may be replaced, up to objectwise weak equivalence, by
a diagram $(\dgrm Z)^{\cat E}$, where \dgrm Z is a fibrant approximation of
the realization of $\hom(V_{\dgrm T}, W)$ as a $D$-diagram, i.e., \dgrm Z is an
$F$-local diagram and therefore the map $\hom(g^{\cat E}, \hom(V_{\dgrm T},
W))$ is a weak equivalence.
\end{proof}
\begin{proposition}\label{conditions:verified}
The class of maps
\[
\Hor(F)=\{\Delta^n\otimes A\otimes \dgrm T \coprod_{\partial\Delta^n\otimes
A\otimes \dgrm T} \partial\Delta^n\otimes B\otimes \dgrm T \rarrow \Delta^n
\otimes B\otimes \dgrm T \; |\; n\geq 0,\; \dgrm T\in \cat O\}
\]
is instrumented with a factorization setup \cal H and a cardinal $\kappa
> |A|+|B|$.
\end{proposition}
\begin{proof}
The factorization setup \cal H is constructed as follows. Any map $u=(u_1,u_2)$
\[
\xymatrix{
   \Delta^n\otimes A\otimes \dgrm T \coprod_{\partial\Delta^n\otimes A\otimes \dgrm T} \partial\Delta^n\otimes B\otimes \dgrm T \ar[r]^<(.5){u_1}\ar[d]      &    \dgrm X \ar[d]  \\
   \Delta^n \otimes B\otimes \dgrm T \ar[r]_{u_2}                                                                                                       &    \dgrm Y         \\
}
\]
is uniquely given by a map of \dgrm T into the nodes of the diagram:
\[
\xymatrix@C=8pt{
            &   &   & \hom(\partial\Delta^n \otimes B, \dgrm X)\ar[dl] \ar[dd] \\
  \hom(\Delta^n \otimes B, \dgrm Y) \ar[dd] \ar[rr] &  & \hom(\partial\Delta^n \otimes B, \dgrm Y)\ar[dd]\\
       &    \hom(\Delta^n \otimes A, \dgrm X) \ar[dl]\ar'[r][rr] &  &\hom(\partial\Delta^n \otimes A, \dgrm X)\ar[dl]  \\
    \hom(\Delta^n \otimes A, \dgrm Y) \ar[rr] &  &\hom(\partial\Delta^n \otimes A, \dgrm Y).  \\
}
\]
Take the `3-dimensional pullback' (the inverse limit) in the diagram above
$\dgrm W_{g,n}$. Then, similarly to the proof of Proposition~
\ref{factorization:setup}, define $H_{g,n}$ to be the set of all maps $u$
above which correspond bijectively, by adjointness, to the set of maps $\cal
O(\dgrm W_{g,n})$; assign $\cal H(g) = \bigcup_{n\geq 0}H_{g,n}$. The
factorization property readily follows.

It is straightforward to check that the domains of maps in $\Hor(F)$ are
$\kappa$-small relative to $\Hor(F)$-cell.
\end{proof}

Proposition~\ref{conditions:verified} verifies for the class $F$ the conditions listed in \ref{conditions}. Hence the results of Section~\ref{loc} apply and we obtain an example $L_F$ of the localization with respect to the \underline{class} of maps $F$.

\appendix
\section{Contractible objects and null homotopies in a model category}\label{null-homotopy}
The purpose of this appendix is to discuss the fundamental notions of
contractible objects and null-homotopic maps in a model category. Using these
notions, we prove Proposition~\ref{no-zero-divisors} in a manner that allows
for immediate generalization to an arbitrary simplicial model category
(satisfying some assumptions), as do the rest of the proofs in
Section~\ref{loc}.

Some authors use these notions in pointed model categories \cite{CDI}, where
the definitions are clear: contractible objects are weakly equivalent to the
zero object and a null map is homotopic to a map which factors through the zero
object. Our aim here is to discuss contractible objects and null maps in any
model category, while generalizing the pointed model category case.
\begin{definition}
A \emph{category with contractible objects} is a pair $(\cat M,\textsl{pt})$,
where \cat M is a model category and \textsl{pt} is a distinguished object in
\cat M, which satisfies that every retract of \textsl{pt} is naturally
isomorphic to \textsl{pt}. This distinguished object is called the
\emph{one-point object} or \emph{singleton}. An object $U$ in \cat M is called
\emph{contractible} if $U$ is weakly equivalent to the one-point object
\textsl{pt}. A map $f\colon A\rarrow X$ is called \emph{null homotopic}
(\emph{nullmap} or just \emph{null}) if it factors up to homotopy (both left
and right) through the one-point object.
\end{definition}

\begin{example}
In any model category \cat M the initial and terminal objects may be chosen
to be one-point objects. They lead, however, to different notions of
contractibility and null maps. If \cat M is a pointed model category
($\emptyset=\ast=0$), then the only object which may be taken as a singleton
is the zero object $\textsl{pt}=0$.

The notions of a singleton, contractible object and nullmap are self-dual. In
the present paper we always take $\textsl{pt}=\ast$, the terminal object in the
category of diagrams of spaces.
\end{example}
\begin{proposition}\label{null}
Let \cat M be a model category with a singleton \textsl{pt} and let $f\colon
A\rarrow X$ be null homotopic.
\begin{enumerate}
\item
If $A$ is a cofibrant object in \cat M, then $f$ factors through a cofibrant
contractible object.
\item
If $X$ is a fibrant object in \cat M, then $f$ factors through a fibrant
contractible object.
\end{enumerate}
\end{proposition}
\begin{proof}
We will prove part (1); the proof of part (2) is dual.

Since $f$ is null homotopic, there exists a map $g=\iota\varphi$, where $A \overset \varphi \longrightarrow \textsl{pt}\overset \iota \longrightarrow X$ homotopic (both left and right) to $f$.  Choose a good cylinder object $A\wedge I$ and a left homotopy $H\colon A\wedge I \rarrow X$ between $f$ and $g$. Then the following solid arrow diagram is commutative,
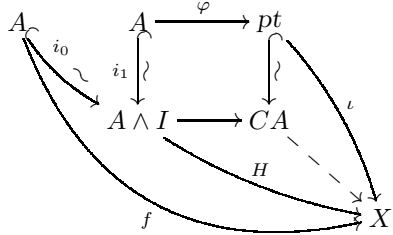
\[
\xymatrix{ A \ar@{^{(}->}@/_/[dr]^<(.2){i_0}^{\dir{~}} \ar@/_3pc/[ddrrr]_f
   &  A \ar[r]^\varphi \ar@{^{(}->}[d]^{\dir{~}}_{i_1} & \textsl{pt}\ar@{^{(}->}[d]^{\dir{~}} \ar@/^/[ddr]^\iota\\
   &  A \wedge I \ar[r] \ar@/_/[drr]^H      & CA \ar@{-->}[dr]\\
   &                                       &         & X}
\]
where $CA=\textsl{pt}\coprod_A A\wedge I$. Then there exists the natural
arrow $CA\rarrow X$ which preserves the diagram commutative.

If $A$ is cofibrant, then the maps $i_0,i_1 \colon  A\rarrow A\wedge I$ are
trivial cofibrations, since $A\wedge I$ is a good cylinder object. Hence, the
cobase change of $i_1$ is also a trivial cofibration, thus the cone object
$CA$ is contractible, and we obtain the required factorization (if $CA$ is not
cofibrant, then factor the map $A\rarrow CA$ as a cofibration followed by an
acyclic fibration and obtain a contractible object $C'A$ through which the map
$f$ factors).

The dual of a cone object is a based paths object.
\end{proof}

\bibliographystyle{abbrv}
\bibliography{xbib}

\end{document}

%% file: symbols.tex
\newcommand{\op}{{\ensuremath{\textup{op}}}}
\newcommand{\Hor}{\ensuremath{\textup{Hor}}}

\newcommand{\dgrm}[1]{\ensuremath{\smash{\underset{\widetilde{\hphantom{#1}}}{#1}} \mathstrut}}

\newcommand{\mhyphen}{\ensuremath{\mathrel{\mathstrut}}}

\newcommand{\domain}[1]{\ensuremath{\mathrm{dom}({#1})}}
\newcommand{\range}[1]{\ensuremath{\mathrm{codom}({#1})}}

\newcommand{\cal}[1]{\ensuremath{\mathcal #1}}

\newtheorem {theorem1}{Theorem}[section]
\newtheorem {theorem}[theorem1]{Theorem}
\newtheorem {corollary}[theorem1]{Corollary}
\newtheorem {proposition}[theorem1]{Proposition}
\newtheorem {lemma}[theorem1]{Lemma}
\theoremstyle{definition}
\newtheorem {definition}[theorem1]{Definition}
\newtheorem {example}[theorem1]{Example}
\theoremstyle{remark}
\newtheorem {remark}[theorem1]{Remark}

\newcommand{\calS}{\ensuremath{\mathcal{S}}}
\newcommand{\calT}{\ensuremath{\mathcal{T}}}

\newcommand{\cat}[1]{\ensuremath{\EuScript #1}}

\newcommand{\Map}[1]{\ensuremath{\textup{Map}\,\EuScript #1}}
\newcommand{\Pow}[1]{\ensuremath{\textup{Pow}\, #1}}
\newcommand{\cofib}{\ensuremath{\hookrightarrow}}

\newcommand{\Top}{\ensuremath{\calT\textup{op}}}

\newcommand{\Sets}{\ensuremath{\cal{S}\textup{ets}}}

\newcommand{\holim}{\ensuremath{\mathop{\textup{holim}}}}
\newcommand{\colim}{\ensuremath{\mathop{\textup{colim}}}}

\def\id{\ensuremath{\textsl{id}}}

\newcommand{\rarrow}{\rightarrow}

\newcommand{\obj}[1]{\ensuremath{\textup{obj}(#1)}}
\newcommand{\mor}[1]{\ensuremath{\textup{mor}(#1)}}
\newcommand{\overcat}{\ensuremath{\mathbin\downarrow}}